\documentclass[11pt]{amsart}
\usepackage{amsmath, amssymb, amscd, mathrsfs, url, pinlabel,setspace,hyperref,verbatim}
\usepackage[margin=1.25in,marginparwidth=1in,centering,letterpaper,dvips]{geometry}
\usepackage{color,dcpic,latexsym,graphicx,epstopdf,comment,todonotes}
\usepackage[all]{xy}
\usepackage{tikz,pgfplots}
\usepackage{tikz-cd,circuitikz}
\usepackage{enumitem,upquote,tabularx,textcomp,setspace,color}
\usepackage{makecell}
\usepackage{caption} 
\xyoption{rotate}

\newtheorem*{rep@theorem}{\rep@title}
\newcommand{\newreptheorem}[2]{%
\newenvironment{rep#1}[1]{%
 \def\rep@title{#2 \ref{##1}}%
 \begin{rep@theorem}}%
 {\end{rep@theorem}}}
\makeatother

\newtheorem {theorem}{Theorem}
\newreptheorem{theorem}{Theorem}
\newtheorem {lemma}[theorem]{Lemma}
\newtheorem {proposition}[theorem]{Proposition}
\newtheorem {corollary}[theorem]{Corollary}
\newtheorem {conjecture}[theorem]{Conjecture}

\newtheorem {question}[theorem]{Question}

\numberwithin{equation}{section}
\numberwithin{theorem}{section}

\theoremstyle{definition}
\newtheorem{definition}[theorem]{Definition}
\newtheorem{construction}[theorem]{Construction}
\newtheorem{data}[theorem]{Data}
\newtheorem{notation}[theorem]{Notation}

\newtheorem{remark}[theorem]{Remark}

\newtheorem{example}[theorem]{Example}

\newtheorem*{ack}{Acknowledgement}
\newtheorem*{org}{Organization}

\newlist{pcases}{enumerate}{1}
\setlist[pcases]{
  label=\bf{Case~\arabic*:}\protect\thiscase.~,
  ref=\arabic*,
  align=left,
  labelsep=0pt,
  leftmargin=0pt,
  labelwidth=0pt,
  parsep=0pt
}
\newcommand{\case}[1][]{%
  \if\relax\detokenize{#1}\relax
    \def\thiscase{}%
  \else
    \def\thiscase{~#1}%
  \fi
  \item
}

\newcommand{\Z}{\mathbb{Z}}

\newcommand{\rk}{\operatorname{rk}}

\newcommand{\R}{\mathbb{R}}
\newcommand{\C}{\mathbb{C}}
\newcommand{\Ct}{\mathbb{C}^\times}

\newcommand{\F}{\mathbb{F}}
\newcommand{\Q}{\mathbb{Q}}

\newcommand{\bbslash}{\backslash\backslash}

\newcommand{\K}{\mathbb{K}}

\newcommand{\ov}[1]{{\overline{#1}}}
\newcommand{\un}[1]{{\underline{#1}}}

\newcommand{\wti}[1]{{\widetilde{#1}}}

\DeclareFontFamily{U}{mathx}{\hyphenchar\font45}
\DeclareFontShape{U}{mathx}{m}{n}{
      <5> <6> <7> <8> <9> <10>
      <10.95> <12> <14.4> <17.28> <20.74> <24.88>
      mathx10
      }{}
\DeclareSymbolFont{mathx}{U}{mathx}{m}{n}
\DeclareFontSubstitution{U}{mathx}{m}{n}
\DeclareMathAccent{\widecheck}{0}{mathx}{"71}

\newcommand{\id}{\operatorname{id}}

\DeclareMathOperator{\im}{im}

\usetikzlibrary{calc,intersections}
\tikzset{every picture/.style=thick}
\tikzset{link/.style = { white, double = black, line width = 1.75pt, double distance = 1.25pt, looseness=1.75 }}
\tikzset{crossing/.style = {draw, circle, dotted, minimum size=0.5cm, inner sep=0, outer sep=0}}
\pgfplotsset{compat=1.12}

\usepackage{extarrows}
\usepackage{mathabx}
\usepackage{amsbsy}
\usepackage{appendix}
\usepackage{float}
\usepackage{subfigure}
\usepackage[numbers,sort&compress]{natbib}
\usepackage{amsthm}
\usepackage{geometry}
\usepackage{fancyhdr}
\usepackage{overpic}
\usepackage{mathrsfs}
\usepackage{colonequals}
\usepackage{microtype}
\usepackage{diagbox}
\usepackage{tabularx}

\newcommand{\bpf}{\begin{proof}}
\newcommand{\epf}{\end{proof}}
\newcommand{\bthm}{\begin{theorem}}
\newcommand{\ethm}{\end{theorem}}
\newcommand{\bprop}{\begin{proposition}}
\newcommand{\eprop}{\end{proposition}}
\newcommand{\bcor}{\begin{corollary}}
\newcommand{\ecor}{\end{corollary}}
\newcommand{\blem}{\begin{lemma}}
\newcommand{\elem}{\end{lemma}}
\newcommand{\bdefn}{\begin{definition}}
\newcommand{\edefn}{\end{definition}}
\newcommand{\bcons}{\begin{construction}}
\newcommand{\econs}{\end{construction}}
\newcommand{\bdata}{\begin{data}}
\newcommand{\edata}{\end{data}}
\newcommand{\bexmp}{\begin{example}}
\newcommand{\eexmp}{\end{example}}
\newcommand{\brem}{\begin{remark}}
\newcommand{\erem}{\end{remark}}
\newcommand{\bnot}{\begin{notation}}
\newcommand{\enot}{\end{notation}}
\newcommand{\benu}{\begin{enumerate}}
\newcommand{\benum}{\begin{enumerate}[leftmargin=*]}
\newcommand{\eenu}{\end{enumerate}}
\newcommand{\beq}{\begin{equation}}
\newcommand{\eeq}{\end{equation}}

\newcommand{\al}{\alpha}
\newcommand{\be}{\beta}
\newcommand{\ga}{\gamma}
\newcommand{\Ga}{\Gamma}

\newcommand{\ot}{\otimes}
\newcommand{\op}{\oplus}

\newcommand{\p}{\prime}
\newcommand{\pp}{{\prime\prime}}

\newcommand{\aand}{~\mathrm{and}~ }

\newcommand{\xra}{\xrightarrow}

\definecolor{lygreen}{HTML}{016646}
\title{Instanton \texorpdfstring{$2$}{2}-torsion and Dehn surgeries}

\author{Zhenkun Li}
\address{Academy of Mathematics and Systems Science\\Chinese Academy of Sciences}
\email{zhenkun@amss.ac.cn}

\author{Fan Ye}
\address{Department of Mathematics\\Harvard University}
\email{fanye@math.harvard.edu}

\makeatletter

\begin{document}

\begin{abstract}
	In our earlier work on $2$-torsion in instanton Floer homology, we considered only integral surgeries on a knot $K\subset S^3$ and showed that the absence of $2$-torsion forces $K$ to be fibered. The present paper extends the result to all rational surgeries. We prove that if the framed instanton homology $I^{\sharp}(S^3_r(K);\mathbb{Z})$ is $2$-torsion-free for some $r\in \mathbb{Q}_+$, then $K$ is an instanton L-space knot and $r>2g(K)-1$. Leveraging this $2$-torsion perspective, we also obtain new small-surgery obstructions: If either $S^{3}_{5}(K)$ or $S^{3}_{11/2}(K)$ is $SU(2)$-abelian, then $K$ must be the unknot or the right-handed trefoil. This result sharpens the small-$SU(2)$-abelian surgery theorems of Kronheimer--Mrowka, Baldwin--Sivek, and Baldwin--Li--Sivek--Ye. 
    
\end{abstract}

\maketitle
\section{Introduction}\label{sec: introduction}
The fundamental group $\pi_1(Y)$ of a $3$-manifold $Y$ is a central topic of study in low-dimensional topology. Rather than analyzing $\pi_1(Y)$ directly, one often probes it via representations into simpler groups. The compact Lie group $SU(2)$ is a natural choice, being one of the simplest non-abelian targets. Any reducible $SU(2)$-representation $\rho:\pi_1(Y)\to SU(2)$ has an abelian image and therefore factors through the abelianization $\mathrm{Ab}(\pi_1(Y))\cong H_1(Y;\Z)$. Consequently, irreducible $SU(2)$-representations detect information that is invisible to $H_1(Y;\Z)$ and capture genuinely non-abelian features of $\pi_1(Y)$.

A closed $3$-manifold $Y$ is \emph{$SU(2)$-abelian} if $\pi_{1}(Y)$ admits no irreducible representation into $SU(2)$. Determining whether a given manifold is $SU(2)$-abelian is therefore a central problem. For instance, the Property P conjecture asserts that the $1$-surgery $S^{3}_{1}(K)$ on a knot $K$ is $SU(2)$-abelian precisely when $K$ is the unknot. More broadly, Problem 3.105(A) of Kirby’s problem list \cite{kirby1997problems} conjectures that an integer homology sphere $Y$ is $SU(2)$-abelian if and only if $Y\cong S^{3}$.

Towards proving the Property P conjecture, Kronheimer and Mrowka \cite{kronheimer04su2} showed that for any rational slope $|r|\leq 2$, the surgery manifold $S^3_r(K)$ is not $SU(2)$-abelian for any non-trivial knot $K$. They also raised the following question.
\begin{question}\label{ques: 3 and 4}
	If $K$ is non-trivial, are $S^3_{3}(K)$ and $S^3_4(K)$ always non-$SU(2)$-abelian?
\end{question}

They only asked for slopes $3$ and $4$ as it is well-known that the $5$-surgery $S^5_5(T_{2,3})$ on the right-handed trefoil $T_{2,3}$ is homeomorphic to the lens space $L(5,4)$ (cf. \cite[Proposition 3.1]{moser1971elementary} for $r=3,s=2,p=1,q=-5$), which has abelian fundamental group and is $SU(2)$-abelian. However, the natural question is whether $T_{2,3}$ is the only exception.

\begin{question}\label{ques: 5}If $K$ is not the unknot or the right-handed trefoil $T_{2,3}$, is $S^3_5(K)$ always non-$SU(2)$-abelian?
\end{question}

Recent developments towards understanding the $\C$-coefficient framed instanton homology $I^\sharp(S^3_r(K);\C)$ of Dehn surgeries on knots enable us to fully answer Question \ref{ques: 3 and 4}: In \cite{baldwin2019lspace}, Baldwin and Sivek resolved the case of $4$-surgery, while in \cite{BLSY21}, the authors of the current paper collaborating with Baldwin and Sivek resolved the case of $3$-surgery.

The idea behind the proofs for the slopes $3$ and $4$ is that if a surgery manifold $Y$ with prime power order $|H_1(Y)|$ is $SU(2)$-abelian, then $Y$ is an instanton L-space (cf. \cite[Remark 1.6]{baldwin2019lspace}), i.e., $\dim I^\sharp(Y;\C)=|H_1(Y)|$. 

Due to \cite{baldwin2019lspace,LY2021large}, if a knot $K$ admits an instanton L-space surgery for a positive slope (usually called an \emph{instanton L-space knot}), then its topology is restricted. In particular, the genus of the knot satisfies $g(K)\leq (r+1)/2$ for any instanton L-space surgery slope $r\in \Q_+$. Thus, when $r= 3$ or $4$, we know that $g(K)\le 2$. The main ingredient of the rest of the proof is then to study instanton L-space knots with genus $1$ or $2$. In addition the results about slopes $3$ and $4$, we have the following theorem from the previous work.

\bthm[{\cite[Corollary 1.5]{BLSY21}}]\label{thm: BLSY results}
For any surgery slope $r\in \frac{1}{2}\Z$ with $|r|<5$, the surgery manifold $S^3_r(K)$ is $SU(2)$-abelian if and only if $K$ is the unknot.
\ethm

From the inequality $g(K)\leq (r+1)/{2}$, the slope $5$ involves the case of genus-$3$ instanton L-space knots, which is much harder to study than the case of genus $1$ or $2$. Thus, all of the old techniques fail when we attempt to tackle Question \ref{ques: 5}.

The study of the $2$-torsion in framed instanton homology $I^\sharp(S^3_r(K);\Z)$ provides new techniques to tackle the problem for slope $r\geq 5$. Some pioneer results regarding the $2$-torsion were developed by the authors of the current paper in \cite{LY2025torsion}, where the authors only worked with integral surgeries, yet in the current paper we extend relevant results to rational surgeries. 

In particular, the main theorem of this paper is the following.

\bthm\label{thm: main}
	Suppose $K\subset S^3$ is a non-trivial knot and $r\in\Q_+$ such that $I^{\sharp}(S^3_r(K);\Z)$ has no $2$-torsion. Then $K$ is an instanton L-space knot, and $r> 2g(K) - 1$. Furthermore, in the following special cases, we have a better lower bound:
	\begin{itemize}
		\item If $r\in\Z_+$, then we have
		\[
			r\geq 2g(K) - 1 + t_2(S^3_1(K)),
		\]
		where for any closed oriented $3$-manifold $Y$, we define
		\begin{equation*}\label{eq: defn of t_2 main}
		    t_2(Y)=\frac{1}{2}\left(\dim I^\sharp(Y;\F_2=\Z/2)-\dim I^\sharp(Y;\C)\right).
		\end{equation*}Note that $t_2(S^3_1(K))\ge 1$ by \cite[Theorem 1.1]{LY2025torsion}.
		\item If $r\in\Z_+-\frac{1}{2}$, then we have \[r> 2g(K).\]
	\end{itemize}
\ethm
\brem\label{rem: mirror knot trick}
Note that $S^3_{-r}(\ov{K})=-S^3_{r}(K)$ and $g(\ov{K})=g(K)$ for the mirror knot $\ov{K}$ of $K$, where the minus sign before the manifold denotes the orientation reversal. Moreover, we have $\dim I^\sharp(-Y;\K)=\dim I^\sharp(Y;\K)$ for any field $\K$, especially for $\K=\F_2$ and $\C$. Hence, a negative surgery $S^3_r(K)$ has no $2$-torsion if and only if the positive surgery $S^3_{-r}(\ov{K})$ has no $2$-torsion, and we can apply Theorem \ref{thm: main} to $\ov{K}$ to obtain the result for a negative slope.
\erem
\brem
	The inequalities in Theorem \ref{thm: main} still seem far from being sharp to the authors. In the future, we plan to return to this problem and obtain a better inequality that would help us to understand which knot admits $SU(2)$-abelian surgery $S^3_r(K)$ for $r\in (5, 7)$.
\erem
A few applications of Theorem \ref{thm: main} follow. We will prove them in \S \ref{subsec: corollaries}.

First, we show that small surgeries of any non-trivial knot always admit $2$-torsion.

\begin{corollary}\label{cor: small surgery}
	Suppose $K\subset S^3$ is a non-trivial knot, and $r\in \Q$ with $|r|\in (0,2g-1]$. Then $I^{\sharp}(S^3_{r}(K);\Z)$ admits $2$-torsion.
\end{corollary}

Next, we extend \cite[Theorem 1.5]{baldwin2019lspace} by imposing a strict inequality instead of the original non-strict one.
\begin{corollary}\label{cor: 2g-1}
	Suppose $K\subset S^3$ is a non-trivial knot. Suppose $r=p/q$ with $p,q>0$ and $\Delta_K(\zeta^2) \neq 0$ for any $p$-th root of unity $\zeta$, where $\Delta_K(t)$ is the Alexander polynomial of $K$ (indeed, when $p$ is a prime power, the condition for roots always holds by \cite[Remark 1.6]{baldwin2016contact}). Then $S^3_r(K)$ is non-$SU(2)$-abelian unless
	\[
	r > 2g(K) - 1.
	\]
	Furthermore, if $r\in\Z_+$, then
	\[
	r \geq 2g(K) - 1 + t_2(S^3_1(K))
	\]
	and if $r\in\Z_+-\frac{1}{2}$, then
	\[
	r> 2g(K).
	\]
\end{corollary}

Corollary \ref{cor: 2g-1} would be enough to fully resolve Question \ref{ques: 5}.
\begin{corollary}\label{cor: 5- and 5.5-surgery}
	Suppose $K\subset S^3$ is a knot and $r\in\frac{1}{2}\Z$ with $|r|<6$ such that $S^3_r(K)$ is $SU(2)$-abelian. Then one of the following cases happens.
	\begin{itemize}
		\item $K$ is the unknot, 
		\item $K$ is the right-handed trefoil and $r\in\{5,\frac{11}{2}\}$.
            \item $K$ is the left-handed trefoil and $r\in \{-5,-\frac{11}{2}\}$.
	\end{itemize}
\end{corollary}

\begin{remark}
    We would like to mention that prior to the first announcement of our current paper, the authors heard that Sudipta Ghosh had an independent argument to show that $S^3_{5}(K)$ is $SU(2)$-abelian except when $K$ is the unknot or the trefoil, though the argument depends on an in-progress work by Aliakbar Daemi, Mike Miller Eismeier, and Christopher Scaduto.
\end{remark}

One could further study the $7$-surgery and obtain the following.
\begin{corollary}\label{cor: 7-surgery}
	Suppose $K\subset S^3$ is a knot such that $S^3_7(K)$ is $SU(2)$-abelian, then either
	\begin{itemize}
		\item $K$ is the unknot or the right-handed trefoil, or
		\item $K$ is a genus-$3$ instanton L-space knot with $t_2(S^3_1(K))\leq 2$.
	\end{itemize}
\end{corollary}

\brem\label{rem: small t_2}
To compare with, in a private discussion, Mike Miller Eismeier computed $I^\sharp(P;\F_2)$ the Poincar\'e sphere $P$, which indicates
\[
t_2(S^3_1(T_{2,3})) = 3.
\]
This indicates that the condition $t_2(S^3_1(K))\leq 2$ might be a tight condition, and to the best of the authors' knowledge, there is no known example of a closed oriented $3$-manifold $Y$ with $t_2(Y)\in \{1,2\}$.
\erem

\begin{org}
    \S \ref{sec: Instanton propositions} is a continuation of the introduction, which describes the strategy of the proof of Theorem \ref{thm: main}. We introduce technical propositions in instanton theory that are  intermediate steps and byproducts of the proof, which may be of independent interest. Other than the main technical result Proposition \ref{prop: 2-torsion free implies being L-space}, the rest propositions in \S \ref{sec: Instanton propositions} will be proven in \S \ref{subsec: illustration of bypass maps} .
    
    \S \ref{sec: first differential} is a preliminary section. We review the construction of first differentials via the contact gluing maps, following the systematic treatment for all slopes (especially for slope $0$) from \cite[\S 6.2]{LY2022integral2}. Instead of simply transporting known results, we also introduce a new way to visualize the instanton knot homology for dual knots in the surgery manifolds, the bypass maps, and the first differentials on the instanton homologies, as discussed in \cite[Remarks 4.15 and 4.38]{LY2020}.

    In \S \ref{sec: tau invariants}, we summarize and improve the results in \cite{ghosh2024tau}, establishing many equivalent definitions of the concordance invariant $\tau_I$, as well as providing a new one that will be used in \S \ref{sec: Proofs of the main results} towards proving our main theorem.

    In \S \ref{sec: Proofs of the main results}, we prove the results presented in \S \ref{sec: Instanton propositions} and the introduction.

    In \S \ref{sec: future direction}, we list some conjectures and open questions for future research.
\end{org}
\begin{ack}
The authors thank Ciprian Manolescu, Jacob Rasmussen, Miller Eismeier, and Christopher Scaduto for helpful comments on the draft of the paper. The authors would like to thank the annual meeting of the Simons Collaboration on New Structures in Low-Dimensional Topology for providing an excellent opportunity and space for discussion and collaboration. The second author is partially supported by Simons Collaboration \#271133 from Peter Kronheimer and is also grateful to Yi Liu for the invitation to BICMR, Peking University while this project is developed. 
\end{ack}
\section{Technical propositions in instanton theory}\label{sec: Instanton propositions}
The strategy to prove Theorem \ref{thm: main} is divided into three parts. The first part deals with the integral case, i.e., when $r=n\in\Z$ and $t_2(S_n^3(K))=0$, we aim to show that $K$ is an instanton L-space knot. The second part reduces the rational case to the integral case by a cabling trick based on \cite[Corollary 7.3]{gordon1983dehn}, i.e., there exists a diffeomorphism.\begin{equation}\label{eq: diffeo for cabling}
	S^3_{pq}(K_{p,q})\cong L(q,p)\#S^3_{p/q}(K),
\end{equation}
where $K_{p,q}$ is the $(p,q)$-cable of $K$, with the convention that $q>1$ denotes the winding number of the cabling pattern. Thus, we can pass the condition of having no $2$-torsion and being an instanton L-space knot back and force between the original knot and its cables. The third part deals with the special cases when $r$ is an integer and a half integer to obtain specific stronger bounds.

In this section, we point out some intermediate steps and byproducts of the first part, which may be of independent interest.

First, we review the definition of dually Floer simple knot and its relation to the $2$-torsion from our previous work \cite{LY2025torsion}.
\bdefn
Suppose $K$ is a knot in $S^3$ and $\wti{K}_{p/q}$ is the dual knot in $S^3_{p/q}(K)$, the surgery manifold of slope $p/q\in \Q$. The knot $K$ is called \emph{dually Floer simple for $n\in\Z\backslash \{0\}$} if \[\dim KHI(S^3_{n}(K),\wti{K}_{n})=\dim I^\sharp(S^3_{n}(K);\C),\]where $KHI$ is the instanton knot homology defined by Kronheimer--Mrowka \cite[Definition 7.13]{kronheimer2010knots} (cf.\ \eqref{eq: KHI}), which is a relative $\Z/2$-graded $\C$-vector space. If we do not specify $n$, then we simply call $K$ dually Floer simple.

Similarly, the knot $K$ is called \emph{rationally dually Floer simple for $p/q\in\Q\backslash\{0\}$} if  \[\dim KHI(S^3_{p/q}(K),\wti{K}_{p/q})=\dim I^\sharp(S^3_{p/q}(K);\C)\] and simply called rationally dually Floer simple if we do not specify $p/q$.
\edefn
\blem[{\cite[Lemma 3.4 and Proposition 3.9]{LY2025torsion}}]\label{lem: 2-tor to C}
If $t_2(S_n^3(K))=0$ for $n\in\Z_{+}$, then $K$ is dually Floer simple for $n+1$.
\elem
\blem[{\cite[Lemma 1.15]{LY2025torsion}}]\label{lem: fibered}
Suppose $K\subset S^3$ is a non-trivial knot. If $K$ is dually Floer simple for $n\in\Z_+$, then $K$ is fibered and V-shaped, and $n>\nu^\sharp(K)=2\tau^\sharp(K)-1>0$, where the notion ``V-shaped" and the concordance invariants $\nu^\sharp$ and $\tau^\sharp$ are defined by Baldwin--Sivek \cite{baldwin2020concordance}.
\elem
By Lemma \ref{lem: 2-tor to C}, we transform the condition about $2$-torsion into the dually Floer simple condition for instanton homologies over $\C$, so that many powerful tools can be applied. The results in Lemma \ref{lem: fibered} are far from being sharp, and actually are all straightforward corollaries of being an instanton L-space knots \cite{baldwin2019lspace,baldwin2020concordance}. 

In \cite[Theorem 1.4 and \S 5]{LY2025torsion}, we proved that an analogous condition of rationally dually Floer simple in Heegaard Floer homology implies that the knot is a Heegaard Floer L-space knot, which provides the original motivation why we conjecture it is also true in instanton theory. However, the proof for Heegaard Floer theory is based on the immersed curve technique that is currently unavailable for instanton homology. This was the main obstruction we were faced with in our previous work. Motivated by the pictures from the immersed curve, we finally manage to find an alternative proof that only depends on sutured theory and thus can be applied to the instanton homology. 

As a result, we obtain the following proposition, which is the main technical result of this paper (and makes Lemma \ref{lem: fibered} a direct corollary). After introducing some notations and assuming some propositions, we will prove it at the end of this section.

\bprop\label{prop: 2-torsion free implies being L-space}
Suppose $K\subset S^3$ is a non-trivial knot. If $K$ is rationally dually Floer simple for $p/q\in\Q_+$, then $K$ is an instanton L-space knot and $p/q>2g(K)-1$.
\eprop
Similar to the proofs in \cite[\S 5]{LY2025torsion} for Heegaard Floer homology, the proof of Proposition \ref{prop: 2-torsion free implies being L-space} depends on the study of the first differentials on the Floer homology of the dual knot. The Heegaard Floer differentials, denoted by $\Psi$ and $\Phi$ in \cite{LY2025torsion}, are obtained from counting holomorphic disks that cover one of the basepoints exactly once in the doubly pointed Heegaard diagram of the knot \cite[\S 4]{sarkar15moving}. The instanton ones are defined in the proof of \cite[Theorem 3.20]{LY2021large} by compositions of bypass maps between instanton knot homologies of various dual knots. The first differentials for the dual knot of slope $-p/q\in\Q\backslash\{0\}$ are given by\begin{equation}\label{eq: d_1 p/q}
    \wti{d}_{1,+}^{p/q},~ \wti{d}_{1,-}^{p/q}:KHI(-S^3_{-p/q}(K),\wti{K}_{-p/q})\to KHI(-S^3_{-p/q}(K),\wti{K}_{-p/q})
\end{equation}\[\mathrm{satisfying}~\left(\wti{d}_{1,\pm}^{p/q}\right)^2=0,\]where $(-S^3_{-p/q}(K),\wti{K}_{-p/q})$ denotes the mirror knot of $\wti{K}_{-p/q}$ and the minus signs are due to technical convention from the contact gluing maps. Since $\wti{K}_{p/q}$ is rationally null-homologous, \cite[Theorem 3.20]{LY2021large} states that the first differentials are on the first pages of the two spectral sequences\begin{equation}\label{eq: spectral sequence}
    KHI(-S^3_{-p/q}(K),\wti{K}_{-p/q})\Longrightarrow I^\sharp(-S^3_{-p/q}(K);\C),
\end{equation}and hence the rationally dually Floer simple condition for $-p/q$ implies that spectral sequences collapse at the first page and thus the first differentials vanish.

In Heegaard Floer setup, the results in \cite[\S 5.2]{LY2025torsion} showed that the local maxima and minima of the immersed curve of the knot complement contribute to the rank of the first differentials, with at most one exception. The vanishing fact of the first differentials implies that the exceptional case must happen, and there are exactly one local maximum and one local minimum, which finally implies that the knot is a Heegaard Floer L-space knot.

In instanton theory, we do not have specific constructions of local maxima and minima. Instead, inspired by the graphical illustration of the bypass maps and the first differentials, we choose to encode the numbers of local maxima and minima by the ranks of the ``first differentials'' for the dual knot of slope $0$:
\begin{equation}\label{eq: d_1 0}
    \wti{d}_{1,+}^{0},~ \wti{d}_{1,-}^{0}:KHI(-S^3_{0}(K),\wti{K}_{0})\to KHI(-S^3_{0}(K),\wti{K}_{0})
\end{equation}\[\mathrm{satisfying}~\left(\wti{d}_{1,\pm}^{0}\right)^2=0.\]Since the $\wti{K}_0$ is not rationally null-homologous, there is no spectral sequence as in \eqref{eq: spectral sequence} for $p/q=0$, but the first differentials are still well-defined because they are just compositions of bypass maps. We will review the definitions of the first differentials for all slopes in \S \ref{subsec: Review of the first differentials}.

Note that from \cite{li2019direct}, \cite[\S 4.2]{LY2020}, and \cite[\S 2.2]{LY2022integral1}, if we fix a Seifert surface of $K$ with minimal genus, then there exists a $\Z$-grading or $(\Z+\frac{1}{2})$-grading on $KHI(S^3_{p/q}(K),\wti{K}_{p/q})$ compatible with the homological $\Z/2$-grading, depending on whether $p$ is odd or even. We call this extra grading by the \emph{Alexander grading} and denote it by $KHI(S^3_{p/q}(K),\wti{K}_{p/q},h)$. The first differentials are homogeneous with respect to the Alexander grading and we write \begin{equation}\label{eq: restriction of d1}
    \wti{d}^{p/q}_{1,\pm}|_h:KHI(-S^3_{-p/q}(K),\wti{K}_{-p/q},h)\to KHI(-S^3_{-p/q}(K),\wti{K}_{-p/q},h\pm \frac{|p|}{2})
\end{equation}for the restrictions of the first differentials.

We have the following proposition that is analogous to \cite[Propositions 5.3 and 5.4]{LY2025torsion}.
\bprop\label{prop: rank d1}
Suppose $p/q\in \Q$ with $p>0$, $h\in \Z$, and $\ov{\tau}=\tau_I(\ov{K})=-\tau_I(K)$. We have the following cases.
\benum
\item If $q>0$, then \[\rk\wti{d}^{p/q}_{1,-}|_{h+\frac{-1+p}{2}}\ge \rk \wti{d}^0_{1,+}|_{h-\frac{1}{2}}-\delta,\]\[\mathrm{where}~\delta=\begin{cases}
    1&\mathrm{if}~h=\ov{\tau}>0,\aand p/q>2\ov{\tau}-1;\\
    0&\mathrm{otherwise}.
\end{cases}\]
\item If $q>0$, then \[\rk\wti{d}^{p/q}_{1,+}|_{h+\frac{1-p}{2}}\ge \rk \wti{d}^0_{1,-}|_{h+\frac{1}{2}}-\delta,\]\[\mathrm{where}~\delta=\begin{cases}
    1&\mathrm{if}~h=-\ov{\tau}<0,\aand p/q>2\ov{\tau}-1;\\
    0&\mathrm{otherwise}.
\end{cases}\]
\item If $q<0$, then \[\rk\wti{d}^{p/q}_{1,+}|_{h+\frac{-1-p}{2}}\ge \rk \wti{d}^0_{1,+}|_{h-\frac{1}{2}}-\delta,\]\[\mathrm{where}~\delta=\begin{cases}
    1&\mathrm{if}~h=-\ov{\tau}>0,\aand p/q<2\ov{\tau}-1;\\
    0&\mathrm{otherwise}.
\end{cases}\]
\item If $q<0$, then \[\rk\wti{d}^{p/q}_{1,-}|_{h+\frac{1+p}{2}}\ge \rk \wti{d}^0_{1,-}|_{h+\frac{1}{2}}-\delta,\]\[\mathrm{where}~\delta=\begin{cases}
    1&\mathrm{if}~h=\ov{\tau}<0,\aand p/q<2\ov{\tau}-1;\\
    0&\mathrm{otherwise}.
\end{cases}\]
\eenu
\eprop
\brem
The invariant $\tau_I$ was defined and studied in \cite{li2019direct,ghosh2024tau}. It is identified with $\tau^\sharp$ by \cite[Theorem 1.2]{ghosh2024tau}, so we will not distinguish the two tau invariants. The reason to use the invariants of the mirror knot $\ov{K}$ instead of the original knot $K$ is that the first differentials are defined on the instanton knot homology of the mirror of the dual knot, and we want to compare the statements with those in \cite[Propositions 5.3 and 5.4]{LY2025torsion}. On the other hand, from \cite[Propositions 5.3 and 5.4]{LY2025torsion}, one may expect that when $\delta=1$, there is some further restriction on the $\epsilon^\sharp$ invariant defined in \cite{baldwin2022concordanceII} that is analogous to the $\epsilon$ invariant in Heegaard Floer theory \cite{Hom2014epsilon}. Since the current restriction is enough for our purpose, we leave this discussion to the future.
\erem

Meanwhile, we have the following proposition for the first differential of slope $0$.
\bprop\label{prop: rank d1 at least}
Suppose $h\in \Z+\frac{1}{2}$. If either $\tau_I(K)\neq 0$ or $h\neq -\frac{1}{2}$, then\[\rk \wti{d}^0_{1,+}|_{h}\ge \frac{1}{2}\left(\dim KHI(-S_0^3(K),\wti{K}_0,h)-\dim KHI(-S_0^3(K),\wti{K}_0,h+ 1)\right).\]If either $\tau_I(K)\neq 0$ or $h\neq \frac{1}{2}$, then \[\rk \wti{d}^0_{1,-}|_{h}\ge \frac{1}{2}\left(\dim KHI(-S_0^3(K),\wti{K}_0,h)-\dim KHI(-S_0^3(K),\wti{K}_0,h- 1)\right).\]
\eprop
\bprop\label{prop: d0 diff}
Suppose $K\subset S^3$ is a knot of genus $g>0$. Then we have \begin{equation*}\label{eq: half rank main}
    \rk \wti{d}^0_{1,\pm}|_{\pm (g-\frac{1}{2})}=\frac{1}{2}\dim KHI(-S_0^3(K),\wti{K}_0,\pm (g-\frac{1}{2}))\ge 1.
\end{equation*}
Hence we have
\begin{equation*}\label{eq: homology zero main}
    H_*(KHI(-S_0^3(K),\wti{K}_0,\pm (g-\frac{1}{2})),\wti{d}^0_{1,\pm}|_{\pm (g-\frac{1}{2})})=0.
\end{equation*}
\eprop
Assuming Propositions \ref{prop: rank d1}, \ref{prop: rank d1 at least}, and \ref{prop: d0 diff}, we prove Proposition \ref{prop: 2-torsion free implies being L-space} as follows.
\bpf[Proof of Proposition \ref{prop: 2-torsion free implies being L-space}]
Suppose $K$ is rationally dually Floer simple for $p/q\in\Q_+$, and $g=g(K)$. The spectral sequences \eqref{eq: spectral sequence} imply that $\rk \wti{d}^{-p/q}_{1,\pm}=0$. Note that here we use $-p/q$ instead of $p/q$.

Meanwhile, from Proposition \ref{prop: d0 diff}, we have \begin{equation}\label{eq: inequalities d 0}
    \rk \wti{d}^0_{1,\pm}\ge \rk \wti{d}^0_{1,\pm}|_{\pm(g-\frac{1}{2})}\ge  1.
\end{equation}By the third and fourth cases in Proposition \ref{prop: rank d1}, the equalities in \eqref{eq: inequalities d 0} must hold, and $\delta=1$. This will imply that $\tau_I(K)=g>0$ and $p/q>2g-1$.

From \cite[Theorem 2.21 (1)]{LY2020}, we have $KHI(-S_0^3,\wti{K}_0,h)=0$ for $|h|>g-\frac{1}{2}$. From Proposition \ref{prop: rank d1 at least} and the equalities in \eqref{eq: inequalities d 0}, we have \[\dim KHI(-S_0^3(K),\wti{K},h)=2\]for $h\in\Z+\frac{1}{2}$ with $|h|\le g-\frac{1}{2}$. Hence we have \[\dim KHI(-S_0^3(K),\wti{K}_0)=4g\aand \mathrm{then}~\dim KHI(S_0^3(K),\wti{K}_0)=4g,\]since the mirror knot induces the dual vector space by \cite[Proposition 3.3]{LY2021large}.

Finally, from \cite[Lemma 2.5 and Remark 2.6]{LY2025torsion} (cf.\ Lemma \ref{lem: V-shaped}) and $\tau_I(K)=g$, we obtain that\[\dim KHI(S^3_{2g}(K),\wti{K}_{2g})=\dim KHI(S^3_0(K),\wti{K}_0)-2g=2g.\]The spectral sequence \eqref{eq: spectral sequence} and the Euler characteristic result\[\chi(I^\sharp(S^3_{2g}(K));\C)=2g\]from \cite[Corollary 1.4]{Scaduto2015instanton} imply that \[\dim I^\sharp(S^3_{2g}(K));\C)=2g.\]Thus, we know that $K$ is an instanton L-space knot.
\epf

\section{First differentials and bypass maps}\label{sec: first differential}

\subsection{Review of the first differentials}\label{subsec: Review of the first differentials}
In this subsection, we fix some notation and review the construction of the first differentials on instanton knot homology.

From now on, we fix a knot $K\subset S^3$ and write $M=S^3\bbslash K=S^3\backslash{\rm int}N(K)$ for the knot complement. Let $-M$ denote the orientation reversal of $M$.

We fix the meridian $\mu$ and the Seifert longitude $\lambda$ on $\partial M$ such that $\mu\cdot \lambda=-1$. If $K$ is oriented, by default we orient $\lambda$ in the same direction as $K$ and orient $\mu$ according to the algebraic intersection number. Given two co-prime integers $p$ and $q$, a simple closed curve on $\partial M$ of slope $p/q$ is the one representing the homology class $p[\mu]+q[\lambda]$. Let $\ga_{p/q}\subset \partial M$ consist of two curves of slope $p/q$ with opposite orientations and let $\Ga_{p/q}=\ga_{-p/q}$.\footnote{This notation $\Gamma_{p/q}$ aligns with our previous notation in \cite{LY2022integral1} and others, while the notation $\gamma_{p/q}$ is a more natural notation regarding the study of dual knots of Dehn surgeries of slope $p/q$.} Let $\wti{K}_{p/q}$ be the dual knot in $S^3_{p/q}(K)$, i.e., the core of the Dehn filling solid torus. Then $\ga_{p/q}$ consists of two meridians of $\wti{K}_{p/q}$. If $(p,q)=(\pm 1,0)$, we write $\ga_\mu$ and $\Ga_\mu$ for the corresponding curves.

Then for any co-prime pair $(p,q)$, we know $(M,\ga_{p/q})$ is a balanced sutured manifold in the sense of \cite[Definition 2.2]{juhasz2006holomorphic} and \cite[Definition 4.1]{kronheimer2010knots}. Let $SHI(M,\ga_{p/q})$ be the sutured instanton homology of $(M,\ga_{p/q})$ from \cite[Definition 7.10]{kronheimer2010knots}. From the definition of the instanton knot homology in \cite[Definition 7.13]{kronheimer2010knots}, we have:\begin{equation}\label{eq: KHI}
    KHI(S^3_{p/q}(K),\wti{K}_{p/q})=SHI(M,\ga_{p/q})\aand KHI(-S^3_{p/q}(K),\wti{K}_{p/q})=SHI(-M,-\ga_{p/q}),
\end{equation}where $(-S^3_{p/q}(K),\wti{K}_{p/q})$ is the mirror knot of $(S^3_{p/q}(K),\wti{K}_{p/q})$. When $M$ is the knot complement, we know that $SHI(M,\ga_{p/q})$ always inherits an Alexander grading, even when $p/q=0$, which we denote by \[SHI(M,\ga_{p/q})=\bigoplus_{h\in \Z+\frac{p-1}{2}}SHI(M,\ga_{p/q},h).\] From \cite[Proposition 3.3]{LY2021large}, there are canonical isomorphisms.\begin{equation}\label{eq: mirror SHI}
    SHI(-M,-\ga,i)\cong SHI(M,-\ga,i)^\vee\cong SHI(M,\ga,-i)^\vee,
\end{equation}where $\vee$ denotes the dual space. Hence the Alexander grading on the mirror knot is opposite to the Alexander grading of the original knot.

\brem\label{rem: unit}
    We will regard $SHI$ as a vector space over $\C$. The naturality problem of its construction was resolved in \cite[\S 9]{BS2015naturality} by introducing the notion of the projectively transitive system. Due to the projectivity, any map between $SHI$ is only defined up to multiplication by a unit in $\C$. We could fix any units for maps, at the cost that any equation of maps only holds up to multiplication by a unit in $\C$. We will point out this ambiguity if it affects the proof.
\erem

The first differentials $d_{1,\pm}$ on instanton knot homology were introduced in \cite[Theorem 3.20]{LY2021large} for a rationally null-homologous knot. We provide an alternative construction as follows, which generalizes to homologically essential knots, especially for $\wti{K}_0\subset S^3_0(K)$. See \cite[\S 6.2]{LY2022integral2} for more details.

We will construct two maps, which we call \emph{(positive and negative) first differentials for slope $p/q$}, i.e., two maps
\begin{equation}\label{eq: d_1}
    \wti{d}_{1,\pm}^{p/q}:KHI(-S^3_{-p/q}(K),\wti{K}_{-p/q})\to KHI(-S^3_{-p/q}(K),\wti{K}_{-p/q}),
\end{equation}
which shift the Alexander grading by $\pm p$ and satisfy the equality \[\left(\wti{d}_{1,\pm}^{p/q}\right)^2=0.\]
\brem\label{rem: minus sign}
From \eqref{eq: mirror SHI}, we can also take the dual spaces and the dual maps to obtain the first differential on $KHI(S^3_{-p/q}(K),\wti{K}_{-p/q})$. Since the Alexander grading is preserved for mirror knots, the dual of the positive differential is still a positive differential.
\erem

From a diffeomorphism $M\cong M\cup_{\partial M}[-1,1]\times T^2$ and the product contact structure on $[-1,0]\times T^2$ with the dividing set $\Ga_{p/q}$ on the two boundary components, we obtain a contact gluing map \[G:SHI(-M,-\Ga_{p/q})\ot SHI(-[0,1]\times T^2,-\Ga_{p/q}\cup -\Ga_{p/q})\to SHI(-M,-\Ga_{p/q}).\]The reason for using manifolds and sutures with opposite orientations originates from \cite{baldwin2016contact} and \cite{honda2009contactInSFH}. Note that
\begin{equation}\label{eq: diffeo}
    (-[0,1]\times T^2,-\Ga_{p/q}\cup -\Ga_{p/q})\cong (-[0,1]\times T^2,-\Ga_{\mu}\cup -\Ga_{\mu}).
\end{equation}
Here we fix a diffeomorphism that sends $p\mu+q\lambda$ to $\mu$ and $p_0\mu+q_0\lambda$ to $\lambda$, where the pair $(p_0,q_0)$ satisfies $q_0p-p_0q=1$ and an additional condition that ensures its uniqueness (cf.\ \cite[Definition 4.2]{LY2020}, where, notationally, the roles of $p$ and $q$ are switched).
\begin{itemize}
    \item If $(p,q)=(\pm 1,0)$, then $(p_0,q_0)=(0,\pm 1)$
    \item If $(p,q)=(0,\pm 1)$, then $(p_0,q_0)=(\mp 1,0)$
    \item If $p,q\neq 0$, then $0\le |p_0|\le |p|,~ 0\le |q_0|\le |q|,~pp_0\le 0$, and $qq_0\le 0$. In particular, if $(p,q)=(n,1)$, then $(p_0,q_0)=(-1,0)$
\end{itemize}
A direct computation shows that the diffeomorphism \eqref{eq: diffeo} sends $\lambda$ to $-p_0\mu+p\lambda$.

The discussion after \cite[Lemma 6.7]{LY2022integral2} showed that \begin{equation}\label{eq: C4}
    SHI(-[0,1]\times T^2,-\Ga_{\mu}\cup -\Ga_{\mu})\cong \C^4.
\end{equation}Moreover, the annuli $A_\mu=[0,1]\times \mu$ and $A_\lambda=[0,1]\times \lambda$ induce a $\Z\oplus \Z$-grading, such that the four dimensions are supported in multi-gradings \begin{equation}\label{eq: multigradings}
    (0,0),(0,0),(0,1),(0,-1),
\end{equation}respectively, up to a global shift. We also consider the annulus $A_{-p_0\mu+p\lambda}=[0,1]\times (-p_0\mu+p\lambda)$, whose relative homology class satisfies\[[A_{-p_0\mu+p\lambda}]=-p_0[A_\mu]+p[A_\lambda].\] From \cite[Theorem 1.12]{li2019decomposition}, the $\Z$-grading associated to $A_{-p_0\mu+p\lambda}$ is a linear combination of the gradings associated to $A_\mu$ and $A_\lambda$, and the four dimensions in \eqref{eq: C4} are supported in gradings $0,0,p,-p$, respectively, up to a global shift.

Furthermore, we know that \eqref{eq: C4} is generated by contact elements. More explicitly, the tight contact structures on $([0,1]\times T^2,\Ga_{\mu}\cup \Ga_{\mu})$ are classified in \cite[Theorem 2.2(3) and \S 5.2]{honda2000classification} and denoted by $\xi_n^\pm$ for $n\in \Z_+$, which are obtained by gluing basic slices. Let $\xi_0$ be the product contact structure. From \cite[Proposition 6.6]{LY2022integral2}, we know that \eqref{eq: C4} is generated by contact elements\begin{equation}\label{eq: contact element}
\theta(\xi_0),\theta(\xi_2^+)=\theta(\xi_2^-),\theta(\xi_1^+),\theta(\xi_1^-).
\end{equation}They are supported in multi-gradings \eqref{eq: multigradings}. From \cite[Proposition 4.13]{LY2021large} and the functoriality of contact elements from \cite[Corollary 1.4]{li2018gluing}, we know $\theta(\xi_n^\pm)=0$ for $n\ge 3$.

\bdefn
Given the slope $p/q$, we assume that $p>0$ or $(p,q)=(0,1)$. Then we define the first differentials $\wti{d}_{1,\pm}^{p/q}$ as follows:\begin{equation}\label{eq: first diff}
    \wti{d}_{1,\pm}^{p/q}=G(-\ot \theta(\xi_1^\pm)),
\end{equation}where we implicitly use \eqref{eq: KHI} and \eqref{eq: diffeo}. 
\edefn
\brem
The diffeomorphisms \eqref{eq: diffeo} associated to $(-p,-q)$ and $(p,q)$ are different, but $(-p)/(-q)=p/q$. That is why we need to fix $p>0$ or $(p,q)=(0,1)$.
\erem
Since the union of two copies of $\xi_1^\pm$ with the same signs provides an overtwisted contact structure, the corresponding contact element vanishes by \cite[Theorem 1.3]{BS2016instanton} and we have $(\wti{d}_{1,\pm}^{p/q})^2=0$ by the functoriality of the contact gluing maps.

\brem
This definition of the first differentials can be easily generalized to $SHI(N,\ga)$ for any other $3$-manifold $N$ with torus boundary and $\ga$ consisting of two parallel curves.
\erem
When $(p,q)\ne (0,1)$, the definition in \eqref{eq: first diff} recovers those in \cite[Theorem 3.20]{LY2021large} that are defined by the compositions of bypass maps. This follows directly from the classification of the tight contact structures on $[0,1]\times T^2$ by Honda \cite{honda2000classification}. We will provide an illustration of this fact in \S \ref{subsec: illustration of bypass maps}.

In particular, we know $\wti{d}_{1,\pm}^{p/q}$ shifts the Alexander grading by $\pm p$, respectively, which is consistent with the gradings of $\theta(\xi_1^\pm)$ associated to $A_{-p_0\mu+p\lambda}$. The spectral sequences in \cite[Theorem 3.20]{LY2021large} imply the dimension inequalities\begin{equation*}\label{eq: dim inequality}
    \begin{aligned}
        \dim H_*\left(KHI(-S^3_{p/q}(K),\wti{K}_{p/q}),\wti{d}_{1,\pm}^{p/q}\right)&=\dim KHI(-S^3_{p/q}(K),\wti{K}_{p/q})-2\rk\wti{d}_{1,\pm}^{p/q}  \\&\ge \dim I^\sharp(-S^3_{p/q}(K)).
    \end{aligned}
\end{equation*}

If $(p,q)=(0,1)$, the spectral sequence does not exist, but the first differentials can still be written as the compositions of two bypass maps. A direct computation shows that they preserve the Alexander grading. We will provide a way to distinguish them in \S \ref{subsec: illustration of bypass maps}.

\subsection{Graphical illustration of bypass maps}\label{subsec: illustration of bypass maps}
In this subsection, we provide an illustrative way to understand the bypass maps and the first differentials. See \cite[Remarks 4.15 and 4.38]{LY2020}. 

Given two slopes $p/q$ and $p^\p/q^\p$ with $q^\p p-p^\p q=1$, there are two bypass maps \[\psi_{\pm,p^\p/q^\p}^{p/q}:SHI(-M,-\Ga_{p/q})\to SHI(-M,-\Ga_{p^\p/q^\p}).\] They are contact gluing maps associated with the basic slices in \cite[\S 3.4]{honda2000classification}. If \[(p/q,p^\p/q^\p)\neq \left(1/(n+1),-1/(-n)\right)\] for some $n\in\Z$, the two bypass maps can be distinguished by the Alexander grading shifts. In particular, when \[(p/q,p^\p/q^\p)=\left(n/1,(-n-1)/(-1)\right),\]the bypass maps $\psi^n_{\pm,n+1}$ shift the Alexander grading by $\mp 1/2$, respectively.\footnote{In this situation, the positive and negative bypass attachments coincide with the positive and negative stabilizations of a Legendrian knot, providing a concrete way to fix the signs.} In general, it is subtle to pin down the signs of the bypass maps since they depend on the convention.

The bypass maps satisfy some exact triangles and commutative diagrams, listed in \cite[Proposition 4.14, Lemmas 4.33 and 4.34]{LY2020}. We provide the following approach to understand these properties more explicitly and systematically.\footnote{The same picture in Heegaard Floer theory comes from the immersed curve invariant perspective.}

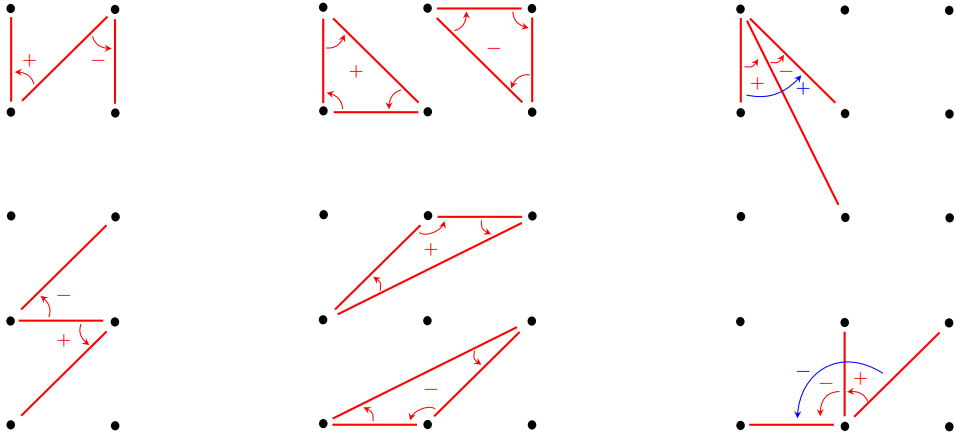
\begin{figure}[ht]
	\centering
	\vspace{-0.5in}
	\begin{overpic}[width = \textwidth]{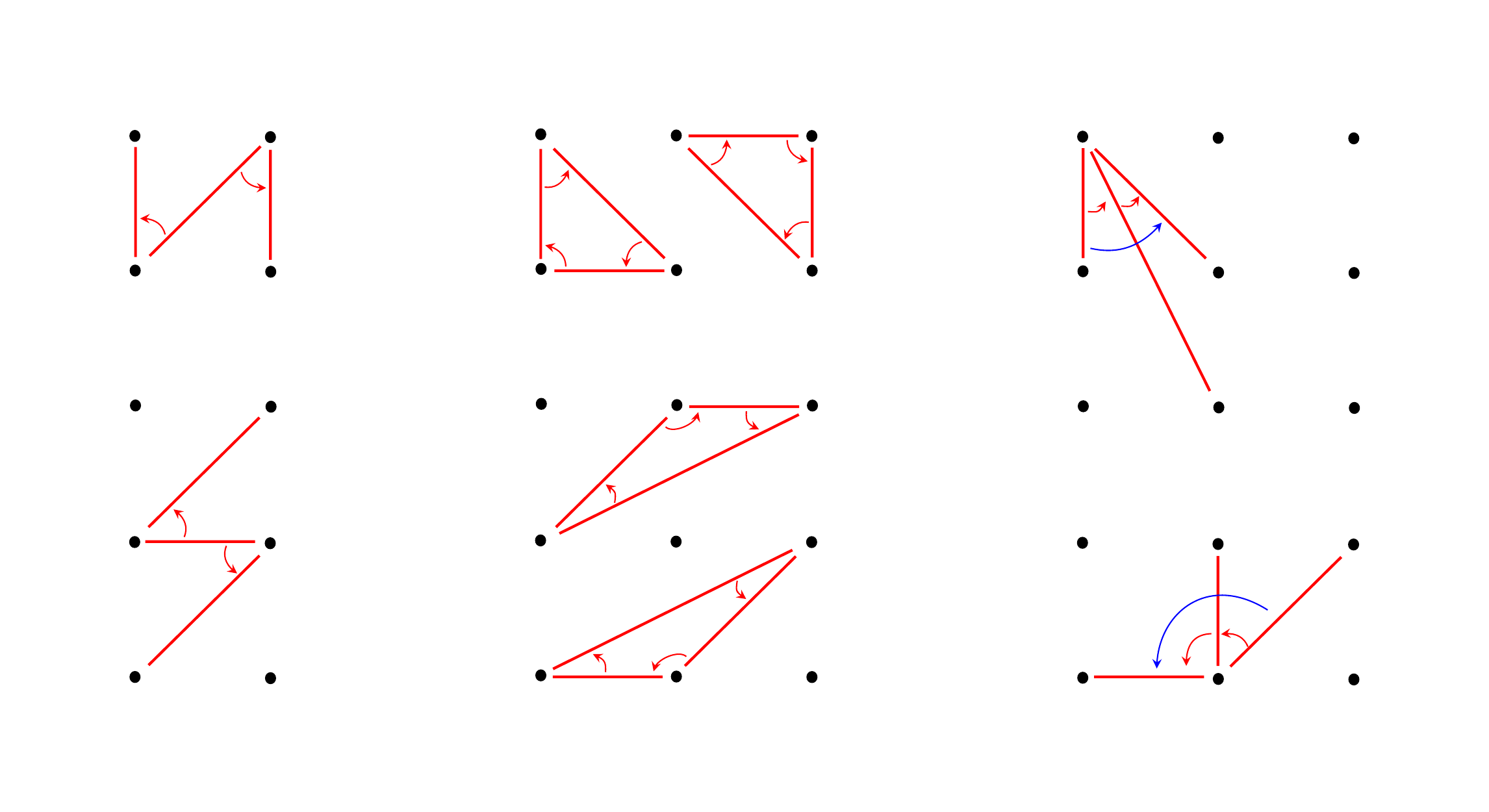}
		\put(10,40.5) {\tiny\color{red}$+$}
		\put(16,40.5) {\tiny\color{red}$-$}
		
		\put(38.5,39.5) {\tiny\color{red}$+$}
		\put(50.5,41.5) {\tiny\color{red}$-$}
		
		\put(73.5,38.5) {\tiny\color{red}$+$}
		\put(76,39.5) {\tiny\color{red}$-$}
		\put(77.5,38) {\tiny\color{blue}$+$}
		
		\put(82.5,12.7) {\tiny\color{red}$+$}
		\put(79.5,12.3) {\tiny\color{red}$-$}
		\put(77.5,13.3) {\tiny\color{blue}$-$}
		
		\put(45,24) {\tiny\color{red}$+$}
		\put(45,11.8) {\tiny\color{red}$-$}
		
		\put(13,16) {\tiny\color{red}$+$}
		\put(13,20) {\tiny\color{red}$-$}
		
	\end{overpic}
	\vspace{-0.5in}
	\caption{Left: two choices of the angle when $h^\p$ has one or two choices; Middle: examples of triangles; Right: examples of compositions. The signs are for bypass maps in \eqref{eq: bypass exact 1} and \eqref{eq: bypass exact 2}.}\label{fig1}
\end{figure}

\benum
    \item\label{t1: arc} Consider the lattice $\Z\times\Z+\frac{1}{2} \subset \R^2$. In the quotient $\R/\Z\times \R$, the lattice becomes a sequence of points, called \emph{basepoints}. We will consider straight arcs connecting basepoints in $\R/\Z\times \R$ and call them \emph{arcs} for short. The second coordinate of the middle point of an arc is called the \emph{height} of the arc. Note that the height of an arc is either in $\Z$ or $\Z+\frac{1}{2}$, depending on whether the numerator $p$ of its slope $p/q\in \Q\cup \{1/0=\infty\}$ is odd or even. It is clearer to draw arcs in $\R^2$ and regard the arcs of the same slope and height as the same. We write $L(p/q,h)$ to denote the arc of slope $p/q$ and height $h$.
    \item\label{t2: bypass} The grading summand of the sutured instanton homology $SHI(-M,-\Ga_{p/q},h)$ corresponds to $L(p/q,h)$ and we will no longer distinguish them. Moreover, we write $L(p/q)$ for $SHI(-M,-\Ga_{p/q})$. For the bypass maps $\psi^{p/q}_{\pm, p^\p/q^\p}$ with $q^\p p-p^\p q=1$, consider two arcs $L(p/q,h)$ and $L(p^\p/q^\p,h^\p)$ that share one endpoint and consider the angle $\al$ in the counterclockwise direction from $L(p/q,h)$ to $L(p^\p/q^\p,h^\p)$. Note that once $h$ is fixed, there are one or two choices of $h^\p$ but always two choices of the angle. See the left subfigure of Figure \ref{fig1}. The two angles correspond to the restrictions of the two bypass maps on the grading summand $h$. Since the signs of the bypass maps are based on convention, it is clearer to denote the bypass map associated to the angle $\al$ by $\psi_{\al}$.
    \item\label{t3: triangle} For two arcs $L(p/q,h)$ and $L(p^\p/q^\p,h^\p)$ that satisfy $q^\p p-p^\p q=1$ and share one endpoint, there is a unique arc connecting the remaining endpoints of the two arcs, denoted by $L(p^\pp/q^\pp,h^\p)$. The assumption $q^\p p-p^\p q=1$ implies that $q^\pp p^\p-p^\pp q^\p=qp^\pp-p q^\pp=1$. Hence there are three bypass maps \[\psi_\al:L(p/q,h)\to L(p^\p/q^\p,h^\p),~\psi_{\al^\p}:L(p^\p/q^\p,h^\p)\to L(p^\pp/q^\pp,h^\p)\]\[\aand \psi_{\al^\pp}:L(p^\pp/q^\pp,h^\p)\to L(p/q,h).\]From \cite[Proposition 4.14]{LY2020}, these three bypass maps form an exact triangle, called the \emph{bypass exact triangle}. In particular, we have the following exact triangles for any $h\in \Z+\frac{n}{2}$ and for any $h\in \Z$, respectively:
\begin{equation}\label{eq: bypass exact 1}
\xymatrix@R6ex{
L(n,h\pm \frac{1}{2})\ar[rr]^{\psi^{n}_{\pm,n+1}}&&L(n+1,h)\ar[dl]^{\psi^{n+1}_{\pm,\mu}}\\
&L(\infty,h\mp \frac{n-1}{2})\ar[ul]^{\psi^{\mu}_{\pm,n}}&
}	
\end{equation} 
\begin{equation}\label{eq: bypass exact 2}
\xymatrix{
L(n-1,h\pm\frac{n}{2})\ar[rr]^{\psi^{n-1}_{\pm,\frac{2n-1}{2}}}&&L(\frac{2n-1}{2},h)\ar[dl]^{\psi^{\frac{2n-1}{2}}_{\pm,n}}\\
&L(n,h\mp\frac{n-1}{2})\ar[ul]^{\psi^{n}_{\pm,n-1}}&
}	
\end{equation}Heuristically, this exact triangle corresponds to the fact that the only basepoints in the triangle formed by $L(p/q,h)$, $L(p^\p/q^\p,h^\p)$, and $L(p^\pp/q^\pp,h^\p)$ are the three vertices and there is \emph{no} basepoints in the interior the triangle. See the middle subfigure of Figure \ref{fig1}. 
    
\begin{figure}[ht]
	\centering
	\vspace{-0.5in}
	\begin{overpic}[width = \textwidth]{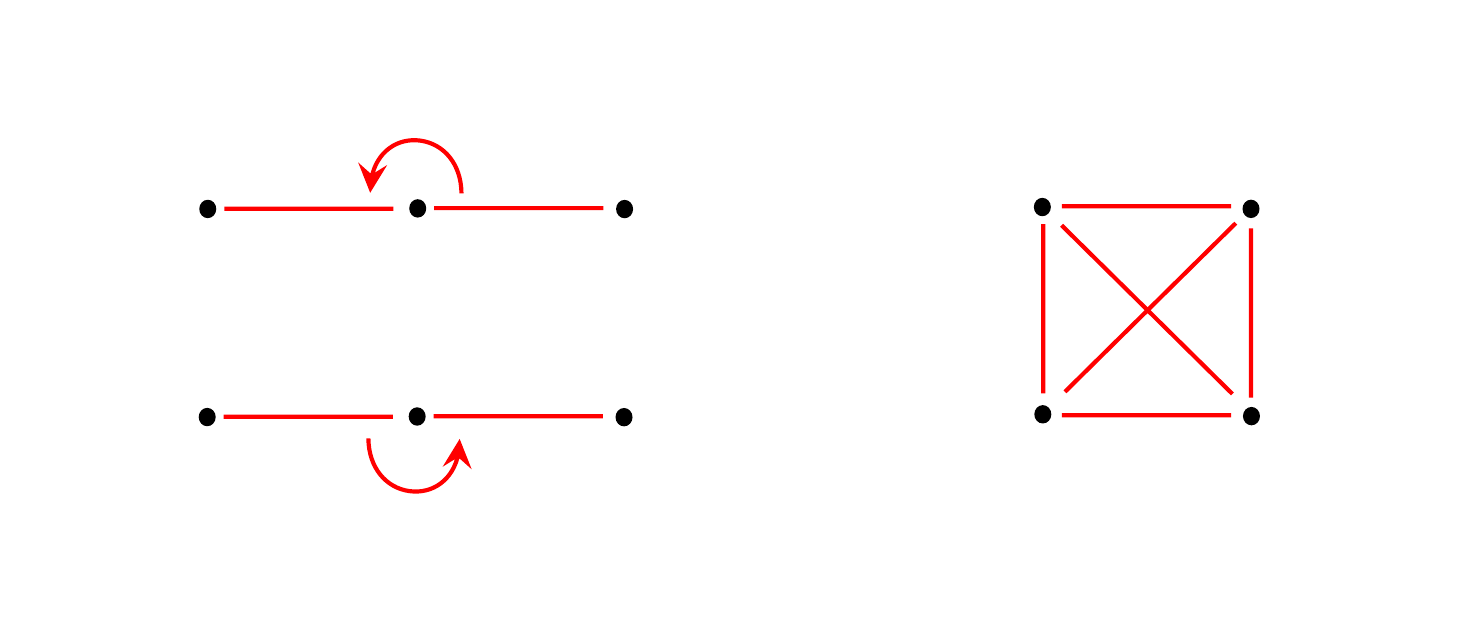}
		\put(27.5,34.5) {\small\color{red}$+$}
		\put(27.5,7.5) {\small\color{red}$-$}
		
		\put(78,12) {\small\color{red}$1$}
		\put(78,30) {\small\color{red}$3$}
		
		\put(69,21) {\small\color{red}$2$}
		\put(87,21) {\small\color{red}$4$}
		
		\put(83.5,17.5) {\small\color{red}$5$}
		\put(81.5,26.0) {\small\color{red}$6$}
	\end{overpic}
	\vspace{-0.6in}
	\caption{Left: the angles corresponding to the first differentials for slope $0$. Right: the six arcs forming a parallelogram, where the numbers denote the index.}\label{fig2}
\end{figure}
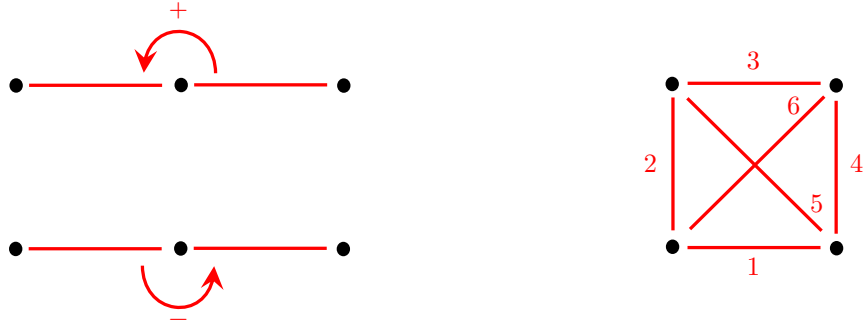

    \item\label{t4: composition}Suppose $L_i$ for $i=1,2,3$ are three arcs that share one endpoint. Suppose the slopes of the arcs satisfy conditions such that there exist bypass maps\[\psi_{\al}: L_1\to L_2,~\psi_{\be}: L_2\to L_3,\aand \psi_{\be\circ \al}: L_1\to L_3\]corresponding to the shared endpoint, where $\be\circ \al$ denotes the composition of the angles. See the right subfigure of Figure \ref{fig1}. From \cite[Lemma 4.34]{LY2020} and Remark \ref{rem: unit}, we have 
    \begin{equation}\label{eq: composition of angle maps}
        \psi_{\be\circ \al}=\psi_{\be}\circ \psi_{\al}
    \end{equation}up to multiplication by a unit. 
    \item\label{t5: angle}Due to \eqref{t4: composition}, for any two arcs $L_1$ and $L_2$ that share one endpoint, the angle $\al$ counterclockwise from $L_1$ to $L_2$ corresponds to a unique map up to multiplication by a unit, obtained from the composition of the bypass maps. We call it a \emph{generalized bypass map} or an \emph{angle map}. Note that a priori, the degree of the angle could be more than $2\pi$. However, from \cite[Proposition 4.13]{LY2021large}, we know such an angle corresponds to the zero map. Thus, we only consider the angle with degree at most $2\pi$. The bypass map is automatically an angle map. The maps \begin{equation}\label{eq: big psi}
        \Psi^{n}_{\pm,n+k}=\psi^{n+k-1}_{\pm,n+k}\circ\cdots \psi^{n}_{\pm,n+1}: SHI(-M,-\Ga_n)\to SHI(-M,-\Ga_{n+k})
    \end{equation}from \cite[Definition 2.10]{LY2022integral1} are also examples of angle maps. We will use $\Psi_\al$ to denote the angle map for the angle $\al$. Note, a priori, one could also mix the signs the sequence of the composition, while in this paper we only work on the case that all signs in the composition are the same.
    \item\label{t6: differential}When $L_1$ and $L_2$ have the same slope, according to the discussion in \S \ref{subsec: Review of the first differentials}, an angle map $\Psi_\al:L_1\to L_2$ could be one of the following cases. 
    \begin{equation}\label{eq: special angle maps}
    	\Psi_\al=\begin{cases}
        \id & \deg \al=0;\\
        \wti{d}^{p/q}_{1,\pm} & \deg \al = \pi;\\
        \wti{d}^{p/q}_{1,+}\circ \wti{d}^{p/q}_{1,-}=\wti{d}^{p/q}_{1,-}\circ \wti{d}^{p/q}_{1,+}&\deg \al=2\pi,
    \end{cases}
    \end{equation}
    where the last equation holds up to multiplication by a unit. We fix the sign convention for the first differentials as follows. For $p/q\neq 0$, the positive differential increases the Alexander grading and the negative differential decreases the Alexander grading. Hence the positive and negative differentials are analogs of $\Psi_h$ and $\Phi_h$ in \cite[Formula (5.7)]{LY2025torsion}, respectively. For $p/q=0$, the positive differential corresponds to the angle above the arc and the negative differential corresponds to the angle below the arc. See the left subfigure of Figure \ref{fig2}.
    \item\label{t7: octahedron} Suppose $L_i$ for $i=1,2,3,4$ are four arcs that form a parallelogram with no basepoint in the interior. In the clockwise direction, the  order of the arcs is $L_1$, $L_2$, $L_3$, and $L_4$. Suppose $L_5$ and $L_6$ are two diagonal arcs of the quadrilateral such that $L_1$, $L_2$ and $L_5$ form a triangle, and $L_2$, $L_3$, and $L_6$ form a triangle. See the right subfigure of Figure \ref{fig2}. Then there exists an octahedral diagram\begin{equation*}
	\xymatrix{
	&L_5\ar[dr]\ar[ddl]&\\
	L_4\ar[ur]\ar[d]&&L_3\ar[ll]\ar[ddl]\\
	L_1\ar[dr]\ar[rr]&&L_2\ar[u]\ar[uul]\\
	&L_6\ar[uul]\ar[ur]&
	}
\end{equation*}where the four exact triangles are from \eqref{t3: triangle}, four commutative triangles are from \eqref{t4: composition}, and the two commutative squares are from \cite[Lemma 4.33]{LY2020}. Note that all commutative diagrams only hold up to multiplication by a unit and ultimately reduce to the classification of the tight contact structure on $[0,1]\times T^2$ from \cite[Lemma 4.21]{honda2000classification} and the interpretation in \cite[\S 1.12]{Honda2000classification2}. Moreover, from diagram chasing or \cite[Proposition 3.6]{LY2022integral1}, we obtain the following two exact triangles\begin{center}
\begin{minipage}{0.4\textwidth}
	\begin{equation}\label{eq: L5 to L6}
\xymatrix{
L_5\ar[rr]&&L_6\ar[dl]\\
&L_2\op L_4\ar[ul]&
}	
\end{equation}
\end{minipage}
\begin{minipage}{0.4\textwidth}
		\begin{equation}\label{eq: L6 to L5}
\xymatrix{
L_6\ar[rr]&&L_5\ar[dl]\\
&L_1\op L_3\ar[ul]&
}	
\end{equation}\end{minipage}
\end{center}where all maps are compositions or direct sum of bypass maps.
\eenu
\brem
The above approach differs from that in \cite[Remarks 4.15 and 4.38]{LY2020} by reflecting $\R^2$ along the second axis, which corresponds to the dual construction in Remark \ref{rem: minus sign}. Indeed, we also need to rotate $\R^2$ by $\pi$, but the fact that $-\ga$ is homotopic to $\ga$ on $\partial M\cong T^2$ implies that the rotation does not affect the illustration. The current approach is more compatible with that in \cite[\S 5]{LY2025torsion} about the immersed curve invariant in Heegaard Floer theory.
\erem

\section{\texorpdfstring{$\tau$}{Tau} invariant}\label{sec: tau invariants}
In this section, we will review many equivalent definitions of the instanton tau invariant $\tau_I$. We then provide a new characterization of tau which we will use in later section to prove our main theorem. We fix a knot $K\subset S^3$ and denote its knot complement by $M=S^3\bbslash K$. We will use the notation for sutures and bypass maps from \S \ref{subsec: Review of the first differentials} without further explanation. We start with the following lemma.



\blem[{\cite[Lemmas 2.16, 2.18 and Proposition 4.1]{LY2022integral1}}]\label{lem: surgery exact triangle}

Suppose $n\in\Z$. There exists an exact triangle\begin{equation}\label{eq: surgery triangle}
	\xymatrix{
	SHI(-M,-\Ga_n)\ar[rr]^{H_n=c_1^n\psi^n_{+,n+1}+c_2^n\psi^n_{-,n+1}}&&SHI(-M,-\Ga_{n+1})\ar[ld]^{F_{n+1}}\\
	&I^\sharp(-S^3)\ar[lu]^{G_n}&
	}
\end{equation}where $c^n_1,c^n_2\in \Ct$ are two constant units depending on $n$. Moreover, there are commutative diagrams up to multiplication by a unit.
\begin{equation}\label{eq: comm F}
	\xymatrix{
	SHI(-M,-\Ga_n)\ar[rd]_{F_n}\ar[rr]^{\psi^n_{\pm,n+1}}&&SHI(-M,-\Ga_{n+1})\ar[ld]^{F_{n+1}}\\
	&I^\sharp(-S^3)&
	}
\end{equation}
\begin{equation}\label{eq: comm G}
	\xymatrix{
	SHI(-M,-\Ga_n)\ar[rr]^{\psi^n_{\pm,n+1}}&&SHI(-M,-\Ga_{n+1})\\
	&I^\sharp(-S^3)\ar[lu]^{G_n}\ar[ru]_{G_{n+1}}&
	}
\end{equation}
\elem
For simplicity, we write $F_n^i$ for the restriction of the map $F_n$ in \eqref{eq: surgery triangle} on grading $i$ of the source, and we write $G_n^i$ for the projection of the map $G_n$ in \eqref{eq: surgery triangle} on grading $i$ of the target.
\bdefn[{\cite{li2019direct,ghosh2024tau}}]\label{defn: tau}
For a knot $K\subset S^3$, we define $\tau_I(K)\in \Z$ by one of the following approaches.
\benum
\item Suppose $\un{\mathrm{KHI}}^-(-S^3,K,i)$ is the direct limit of \[\cdots \to SHI(-M,-\Ga_n)\xra{\psi^n_{-,n+1}} SHI(-M,-\Ga_{n+1})\xra{\psi^{n+1}_{-,n+2}}SHI(-M,-\Ga_{n+2})\to \cdots,\] and let $U$ be defined by the collection of maps $\{\psi^{n}_{+,n+1}\}_{n\in \Z}$, using the fact that \[\psi^{n+1}_{-,n+2}\circ \psi^{n}_{+,n+1}=\psi^{n+1}_{+,n+2}\circ \psi^{n}_{-,n+1}\] from \cite[Lemma 4.33]{LY2020}. Define \[\tau_I(K)=\max\{i\mid\exists x\in \un{\mathrm{KHI}}^-(-S^3,K,i) \text{  such that for any }j\ge 0, ~ U^jx\neq 0\}.\] 
\item We define $\tau_I(K)$ by\[\begin{aligned}
    -2\tau_I(K)=&\max\{n\mid G_n\neq 0\}+1=\min\{n\mid G_n=0\}\\=&\min\{n\mid F_n\neq 0\}-1=\max\{n\mid F_n=0\}
\end{aligned}\]
\item Suppose $n\in \Z$ such that $F_n\neq 0$. Define\[\tau_I(K)=\max\{i\mid F_n^i\neq 0\}-\frac{n-1}{2}=-\min\{i\mid F_n^i\neq 0\} + \frac{n-1}{2}.\]
\item Suppose $n\in \Z$ such that $G_{n}\neq 0$. Define\[\tau_I(K)=\max\{i\mid G_{n}^i\neq 0\}-\frac{n+1}{2}=-\min\{i\mid G_{n}^i\neq 0\} + \frac{n+1}{2}.\]
\eenu
\edefn
\blem[{\cite[Lemma 2.5 and Remark 2.6]{LY2025torsion}}]\label{lem: V-shaped}
The sequence \[\big(\dim KHI(S^3_n(K),\wti{K}_n)\big)_{n\in\Z}\] is V-shaped with unique minimum $2\tau_I(K)$. The sequence \[\big(\dim SHI(-M,-\Ga_n)\big)_{n\in\Z}\] is V-shaped with unique minimum $-2\tau_I(K)$. Note that the minus sign comes from the fact that $\Ga_n=\ga_{-n}$.
\elem
\blem\label{lem: equivalent of tau}
The approaches to define $\tau_I(K)$ in Definition \ref{defn: tau} are equivalent and independent of the choice of $n$.
\elem
\bpf
Indeed, the first definition is the original one and the rest are properties of $\tau_I(K)$. The second definition follows from Lemmas \ref{lem: surgery exact triangle} and \ref{lem: V-shaped}. Note that $\dim I^\sharp (-S^3)=1$.



From \cite[Proposition 5.3]{ghosh2024tau}, we know the third definition holds for sufficiently large $n$. The last definition also holds for sufficiently large $n$ due to the following commutative diagram\begin{equation}\label{eq: dual commutative diagram 2}
    \xymatrix{SHI(-S^3\bbslash K,-\Ga_{n}(K))\ar[rr]^{\cong}&&SHI(-S^3\bbslash \ov{K},-\Ga_{-n}(\ov{K}))^\vee\\
    \\I^\sharp(-S^3)\ar[uu]^{G_{n}(K)}\ar[rr]^\cong &&I^\sharp(-S^3)^\vee\ar[uu]_{F_{-n}(\ov{K})^\vee}}
\end{equation}where the horizontal maps are canonical isomorphisms from \eqref{eq: mirror SHI}. To prove it also holds for small $n$, we introduce the following notation as in \cite[\S 5]{ghosh2024tau}
\[l_n=\begin{cases}
   \max\{i\mid G_{n}^i\neq 0\}-\min\{i\mid G_{n}^i\neq 0\}+1&\mathrm{if}~G_n\neq 0;\\0&\mathrm{if}~G_n=0.
\end{cases}\] From \cite[Lemma 5.9 and Corollary 4.7]{ghosh2024tau}, for sufficiently large $n$, we have \begin{equation*}\label{eq: l}
    l_{-n}=2\tau_I(\ov{K})+n=-2\tau_I(K)+n.
\end{equation*}From \cite[Lemma 5.10]{ghosh2024tau}, if $l_n>0$, then $l_{n+1}\le l_n-1$ and if $l_n=0$, then $l_{n+1}=0$. Then we have \begin{equation}\label{eq: sequence inequality}
    l_{-2\tau_I(K)-1}\le l_{-2\tau_I(K)-2}-1\le \cdots\le l_{-2\tau_I(K)-(-2\tau_I(K)+n)}-(-2\tau_I(K)+n-1)=1.
\end{equation}The second definition implies $G_{-2\tau_I(K)-1}\neq 0$ and then $l_{-2\tau_I(K)-1}\ge 1$. Hence all inequalities in \eqref{eq: sequence inequality} are indeed equalities. From the proof of \cite[Proposition 5.3]{ghosh2024tau}, we know the maximum and the minimum are opposite numbers. Note that here the gradings are in $\Z$ or $\Z+\frac{1}{2}$. Thus, the fourth definition holds for any $n\in\Z$. From \eqref{eq: dual commutative diagram 2}, the third definition holds for any $n\in\Z$.
\epf

We will mainly use the following characterization of $\tau_I(K)$.

\bcor\label{cor: char of tau G}
Suppose $n\in\Z$ and $x$ is a nonzero homogeneous element in $SHI(-M,-\Ga_n,i_0)$. Then the following conditions are equivalent.
\benum
\item $n=-2\tau_I(K)-1$ and $x$ is the generator of $\im G^0_n$. In such a case, we have $i_0=0$.
\item $x$ is a generator of the $1$-dimensional space $\im G_n$.
\item $\psi^n_{\pm,n+1}(x)=0$
\item $x\in \im \psi^\mu_{\pm,n}$.
\eenu
\ecor
\bpf
The equivalence $(1)\Leftrightarrow(2)$ follows from the fourth definition of $\tau_I(K)$ in Definition \ref{defn: tau}. The equivalence $(3)\Leftrightarrow (4)$ follows from the bypass exact triangle in \eqref{eq: bypass exact 1}. The result $(1)\Rightarrow (3)$ follows from the second definition of $\tau_I(K)$ and \eqref{eq: comm G}. The result $(3)\Rightarrow (2)$ is from the surgery exact triangle \eqref{eq: surgery triangle}. 
\epf
\brem
Since the maps $F_n$ and $G_n$ are dual in the sense of \eqref{eq: dual commutative diagram 2}, one may expect a suitable dual result of Corollary \ref{cor: char of tau G} about $F_n$ to hold. However, the authors notice that the statement and the proof for $F_n$ are much more subtle than those for $G_n$, with the key difference that $\im G_n$ is a concrete $1$-dimensional subspace of $SHI(-M,-\Ga_n)$ while we only have a quotient $1$-dimensional space $SHI(-M,-\Ga_n)/ \ker F_n$. This is why we consider only the characterization of $\tau_I$ with respect to $G_n$ and avoid the discussion involving $F_n$.
\erem
\section{Proofs of the main results}\label{sec: Proofs of the main results}
\subsection{Propositions}
In this subsection, we prove the propositions in \S \ref{sec: Instanton propositions}.

Comparing the illustration in \S \ref{subsec: illustration of bypass maps} and that in \cite[\S 5]{LY2025torsion}, the main difference is that currently we use the arc $L$ itself to represent the pairing with some immersed curve obtained from the knot complement. This is because the immersed curve theory has not been developed in instanton theory. In particular, we cannot study the immersed curve by its components and we need to find alternative statements and proofs for results in \cite[\S 5.2]{LY2025torsion}.

The results in \cite[\S 5.2]{LY2025torsion} are about local maxima and local minima of the immersed curve. Due to the description of $\wti{d}^{0}_{1,\pm}$ in \eqref{t6: differential}, the positive differential of slope $0$ captures the local maxima and the negative differential of slope $0$ captures the local minima. More explicitly, the rank of $\wti{d}^{0}_{1,\pm}$ in grading $h$ is the analog of the number $n_\pm (\ga,h\pm 1/2)$ of local maxima or minima.

With the above observation, we will first establish Proposition \ref{prop: rank d1}. We will prove only the first case. Then one can apply the first case to the orientation reversal of the original knot, the mirror knot, and the orientation reversal of the mirror knot to obtain the remaining cases. We restate the proposition for reader's convenience.

\noindent {\bf Proposition \ref{prop: rank d1}, item (1)}
{\it
Suppose $p,q>0$ are co-prime numbers, $h\in \Z$, and $\ov{\tau}=\tau_I(\ov{K})=-\tau_I(K)$. Then \[\rk\wti{d}^{p/q}_{1,-}|_{h+\frac{-1+p}{2}}\ge \rk \wti{d}^0_{1,+}|_{h-\frac{1}{2}}-\delta,\]\[\mathrm{where}~\delta=\begin{cases}
    1&\mathrm{if}~h=\ov{\tau}>0,\aand p/q>2\ov{\tau}-1;\\
    0&\mathrm{otherwise}.
\end{cases}\]
}


\bpf
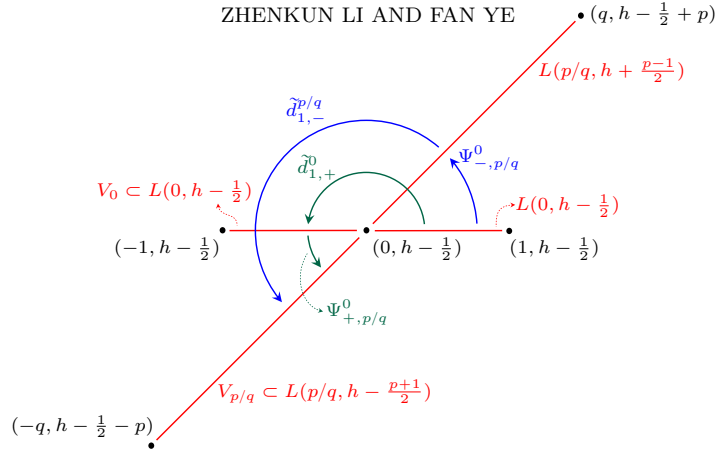
\begin{figure}[ht]
\vspace{-0.6in}
	\centering
	\begin{overpic}[width = 0.5\textwidth]{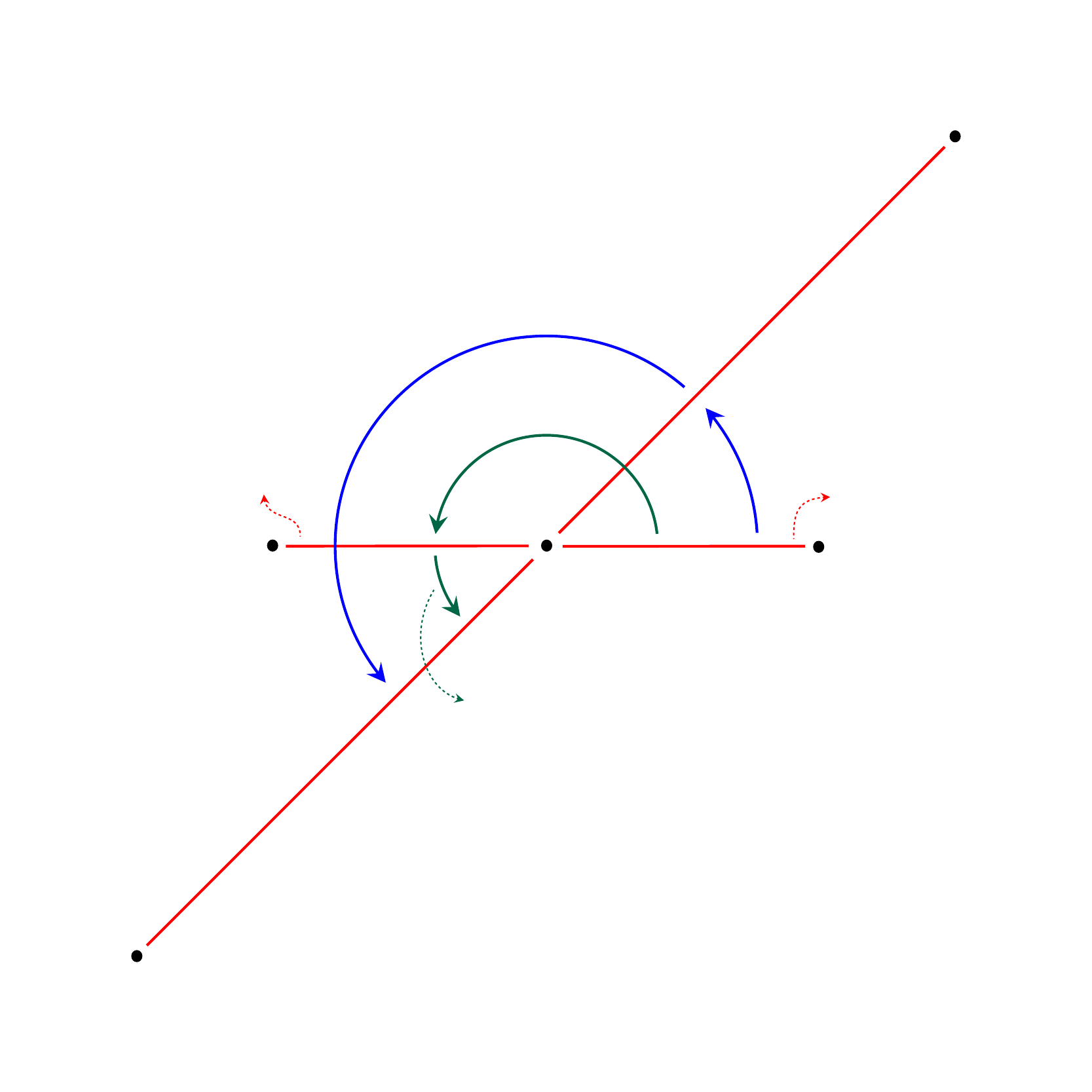}
		\put(89,87) {\tiny $(q,h-\frac{1}{2}+p)$}
		\put(51,46) {\tiny $(0,h-\frac{1}{2})$}
		\put(75,46) {\tiny $(1,h-\frac{1}{2})$}
		\put(6,46) {\tiny $(-1,h-\frac{1}{2})$}
		\put(-12,15) {\tiny $(-q,h-\frac{1}{2}-p)$}
		
		\put(36,70) {\tiny\color{blue}$\wti{d}^{p/q}_{1,-}$}
		\put(66,62) {\tiny\color{blue}$\Psi^{0}_{-, p/q}$}
		\put(38,60) {\tiny\color{lygreen}$\wti{d}^0_{1,+}$}
		\put(43,35) {\tiny\color{lygreen}$\Psi^{0}_{+,p/q}$}
	
		\put(24,21) {\tiny\color{red}$V_{p/q}\subset L(p/q,h-\frac{p+1}{2})$}
		\put(80,77) {\tiny \color{red} $L(p/q, h+\frac{p-1}{2})$}
		\put(3,56) {\tiny\color{red} $V_0\subset L(0,h-\frac{1}{2})$}
		\put(76,53.5) {\tiny\color{red} $L(0,h-\frac{1}{2})$}
		
	\end{overpic}
	\vspace{-0.4in}
	\caption{The angle maps related to $L(p/q)$ and $L(0)$.}\label{fig3}
\end{figure}

For simplicity, we write $L(p/q,h)$ for $SHI(-M,-\Ga_{p/q},h)$, $L(p/q)$ for $SHI(-M,-\Ga_{p/q})$, and adopt the notation of angle maps as in \S \ref{subsec: illustration of bypass maps}. In particular, the maps $\wti{d}^{p/q}_{1,-}|_{h+\frac{-1+p}{2}}$ and $\wti{d}^0_{1,+}|_{h-\frac{1}{2}}$ are the angle maps (with angle $\pi$) as in Figure \ref{fig3}. Also, as in the figure, there are two other angle maps
\[
\Psi^{0}_{-,p/q}: L(0,h-\frac{1}{2})\to L(p/q, h + \frac{p-1}{2}) \text{ and } \Psi^{0}_{+,p/q}: L(0,h-\frac{1}{2})\to L(p/q, h - \frac{p+1}{2}).
\]
We define
\[V_0=\im \wti{d}^0_{1,+}|_{h-\frac{1}{2}}\subset L(0,h-\frac{1}{2})\text{ and } V_{p/q} = \Psi^0_{+,p/q}(V_0).\]
As in Figure \ref{fig3}, there is a commutative diagram according to \eqref{eq: composition of angle maps} (up to multiplication by a unit)
\[
\wti{d}^{p/q}_{1,-}|_{h+\frac{-1+p}{2}}\circ \Psi^0_{-,p/q} = \Psi^0_{+,p/q}\circ \wti{d}^0_{1,+}|_{h-\frac{1}{2}}.
\]
Thus, we have
\[
\begin{aligned}
	V_{p,q} &= \Psi^0_{+,p/q}(V_0)\\
	&= \Psi^0_{+,p/q}\circ \wti{d}^0_{1,+}|_{h-\frac{1}{2}}\left(L(0,h-\frac{1}{2})\right)\\
	&= \wti{d}^{p/q}_{1,-}|_{h+\frac{-1+p}{2}}\circ \Psi^0_{-,p/q}\left(L(0,h-\frac{1}{2})\right),
\end{aligned}
\]
which implies
\[V_{p,q}\subset \im\wti{d}^{p,q}_{1,-}|_{h+\frac{-1+p}{2}}\text{, since } \Psi^0_{-,p/q}\left(L(0,h-\frac{1}{2})\right)\subset L(p/q,h+\frac{-1+p}{2}).\]
As a consequence, the proof of the proposition further reduces to comparing the dimension of the two spaces $V_{p/q}$ and $V_{0}$. The proof is then divided into two steps:
\benum
	\item We first handle the case $q = 1$, i.e., the case when $n=p/q\in\mathbb{Z}$. In this case, since $V_n = \Psi^{0}_{+,n}(V_0)$, we will study the kernel of the map $\Psi^{0}_{+,n}$ directly.
	\item We then handle the case $p/q$ is not an integer. For this case, we find an integer $k$ such that $k<p/q<k+1$ and we will show that
	\[
		\dim V_{k+1} \leq \dim V_{p/q} \leq \dim V_{k},
	\]
	which concludes the proof of the proposition.
\eenu

We will proceed with the case $n=p/q\in\mathbb{Z}$ first. Note that, in this case, we have
\[
\begin{aligned}
	\rk \wti{d}^n_{1,-}|_{h+\frac{-1+n}{2}} &= \dim \im\wti{d}^n_{1,-}|_{h+\frac{-1+n}{2}}\\
	&\geq \dim V_n\\
	&= \dim V_0 - \dim{\rm ker}(\Psi^{0}_{+,n}|_{V_n})\\
	&= \rk \wti{d}^0_{1,+}|_{h-\frac{1}{2}} - \dim{\rm ker}(\Psi^{0}_{+,n}|_{V_n})
\end{aligned}
\]
Thus, it suffices to show that 
\begin{equation}\label{eq: dim ker Psi 0,+,n}
	\dim{\rm ker}(\Psi^{0}_{+,n}|_{V_n}) = \begin{cases}
	1 & \mathrm{if}~ n > 2\ov{\tau} -1;\\
	0 & \mathrm{if}~ n\leq 2\ov{\tau} -1.
\end{cases}
\end{equation}
Recall that, by the definition of the map $\Psi^0_{+,n}$ in \eqref{eq: big psi}, we have
\[
\Psi^0_{+,n} = \psi^{n-1}_{+,n}\circ\cdots\circ \psi^0_{+,1}.
\]
Thus, we can inductively define $V_k = \psi^{k-1}_{+,k}(V_{k-1})$ for $k=1,2,\dots,n-1$. By definition,
\[
\begin{aligned}
	V_0 &= \im \wti{d}^0_{1,+}|_{h-\frac{1}{2}}\\
	&= \psi^{\mu}_{-,0}(\im \psi^{0}_{-,\mu}|_{h-\frac{1}{2}})\\
	&\subset \im \psi^{\mu}_{-,0}.
\end{aligned}
\]

\begin{figure}[ht]
\vspace{-0.2in}
	\centering
	\begin{overpic}[width = 0.9\textwidth]{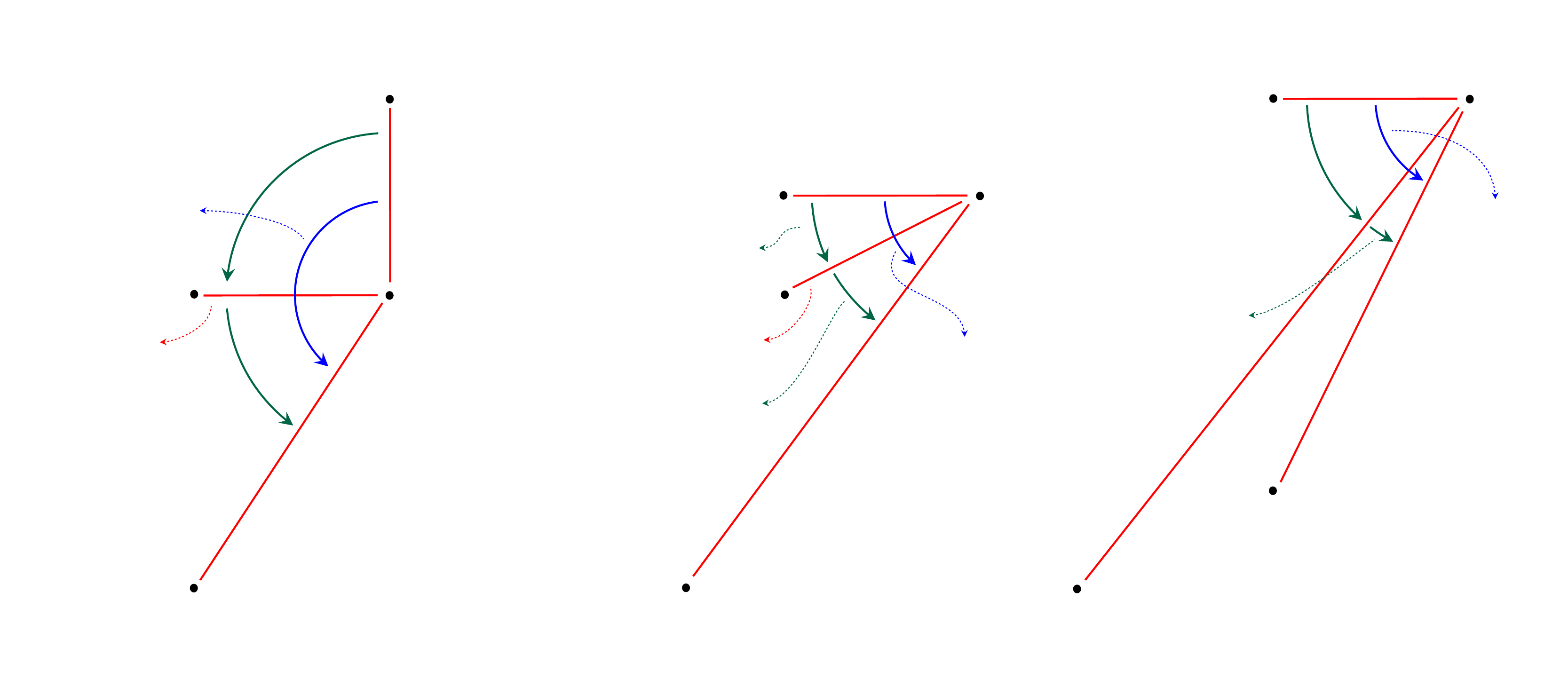}
		\put(25.5,24.5) {\tiny $(0,h-\frac{1}{2})$}
		\put(1,24.5) {\tiny $(-1,h-\frac{1}{2})$}
		\put(25.5,37.5) {\tiny $(0,h+\frac{1}{2})$}
		\put(1,4) {\tiny $(-1,h-\frac{1}{2}-k)$}
		
		\put(25.5,33) {\tiny\color{red} $L(\infty,h)$}
		\put(0,21.5) {\tiny\color{red} $L(0,h-\frac{1}{2})$}
		\put(15,8) {\tiny\color{red} $L(k,h-\frac{k+1}{2})$}
		
		\put(8,30) {\tiny\color{blue} $\psi^{\mu}_{-,n}$}
		\put(16,35) {\tiny\color{lygreen} $\psi^{\mu}_{-,0}$}
		\put(12,18) {\tiny\color{lygreen} $\Psi^{0}_{+,k}$}
		
		\put(39,31) {\tiny $(-1,h-\frac{1}{2})$}
		\put(63,31) {\tiny $(0,h-\frac{1}{2})$}
		\put(37,24.5) {\tiny $(1,h-\frac{1}{2}-k)$}
		\put(38,4) {\tiny $(-q, h-\frac{1}{2}-p)$}
		
		\put(51,32.5) {\tiny\color{red} $L(0,h-\frac{1}{2})$}
		\put(35.5,21.5) {\tiny\color{red} $L(k,h-\frac{k+1}{2})$}
		\put(47.5,10) {\tiny\color{red} $L(p/q,h-\frac{p+1}{2})$}
		
		\put(44,17) {\tiny\color{lygreen} $\Psi^k_{+,p/q}$}
		\put(44,27.5) {\tiny\color{lygreen} $\Psi^0_{+,k}$}
		\put(58,21) {\tiny\color{blue} $\Psi^0_{+,p/q}$}
		
		\put(94.5,37) {\tiny $(0,h-\frac{1}{2})$}
		\put(70,37) {\tiny $(-1,h-\frac{1}{2})$}
		\put(75,10) {\tiny $(-1, h-\frac{3}{2}-k)$}
		\put(60,4) {\tiny $(-q, h-\frac{1}{2}-p)$}
		
		\put(82,38.5) {\tiny\color{red} $L(0,h-\frac{1}{2})$}
		\put(85,18) {\tiny\color{red} $L(k+1,h-\frac{k+3}{2})$}
		\put(61,15) {\tiny\color{red} $L(p/q,h-\frac{p+1}{2})$}
		
		\put(92,30) {\tiny\color{blue} $\Psi^0_{+,k+1}$}
		\put(74,22.5) {\tiny\color{lygreen} $\Psi^{p/q}_{+,k+1}$}
		\put(78,32.5) {\tiny\color{lygreen} $\Psi^{0}_{+,p/q}$}
		
	\end{overpic}
	\vspace{-0.2in}
	\caption{Left: The angle maps related to $L(\infty)$, $L(0)$, and $L(k)$. Middle: the angle maps related to $L(0)$, $L(k)$, and $L(p/q)$. Right: the angle maps related to $L(0)$, $L(p/q)$, and $L(k+1)$.}\label{fig4}
\end{figure}
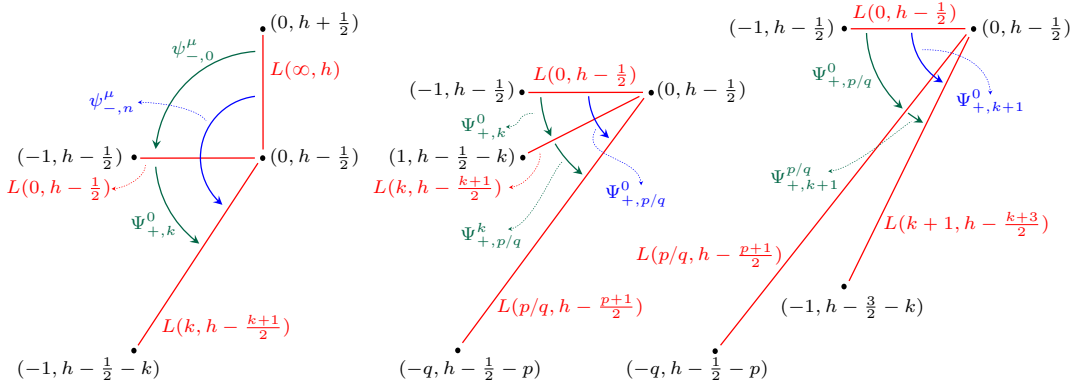

Note that we have an identity $\psi^{\mu}_{-,k} = \Psi^{0}_{+,k}\circ \psi^{\mu}_{-,0}$ as a special case of \eqref{eq: composition of angle maps} that is illustrated in the left subfigure of Figure \ref{fig4}. Thus, we conclude inductively that
\[
V_{k}\subset \im \psi^{\mu}_{-,k} \text{ for all } k\in [0,n].
\]
As a consequence, the exact triangle \eqref{eq: bypass exact 1} implies that $V_k\subset \ker\psi^{k}_{-,k+1}$ for all $k$. Then we conclude that
\[
\ker(\psi^k_{k+1}|_{V_k}) \subset \ker\psi^k_{+,k+1}\cap\ker\psi^k_{-,k+1}.
\]
Then Corollary \ref{cor: char of tau G} implies that
\[
\ker(\psi^k_{k+1}|_{V_k}) = \begin{cases}
	1 & \mathrm{if}~k = 2\ov{\tau} - 1; \\
	0 & \mathrm{otherwise}.
\end{cases}
\]
and Equation \eqref{eq: dim ker Psi 0,+,n} follows. This completes the proof of the case when the surgery slope is an integer.

Next, we consider the general case where $p/q>0$ but $p/q\notin \Z$. Take $k=\lfloor p/q\rfloor$. We have $k<p/q<k+1$. The middle subfigure of Figure \ref{fig4} implies that
\[
\begin{aligned}
	V_{p/q} &= \Psi^0_{+,p/q}(V_0)\\
	&= \Psi^{k}_{+,p/q}\circ \Psi^{0}_{+,k}(V_0)\\
	&= \Psi^k_{+,p/q}(V_k)
\end{aligned}
\]
Thus, we conclude that $\dim V_{p/q}\leq \dim V_k$. Similarly, using the right subfigure of Figure \ref{fig4}, we have
\[
\begin{aligned}
	V_{k+1} &= \Psi^0_{+,k+1}(V_0)\\
	&= \Psi^{p/q}_{+,k+1}\circ \Psi^0_{+,p/q}(V_0)\\
	&=\Psi^{p/q}_{+,k+1}(V_{p/q}). 
\end{aligned}
\]
Thus, we conclude that $\dim V_{k+1}\leq \dim V_{p/q}$. Combined with the above argument, we obtain the desired inequalities
\[
	\dim V_{k+1} \leq \dim V_{p/q} \leq \dim V_{k}.
\]
Thus, the case when $p/q$ is non-integral follows from the results for integers $k$ and $k+1$, and we conclude the proof of the proposition.
\epf
Then we prove the remaining propositions in \S \ref{sec: Instanton propositions}.

\vspace{0.15in}
\noindent{\bf Proposition \ref{prop: rank d1 at least}}
{\it
Suppose $h\in \Z+\frac{1}{2}$. If either $\tau_I(K)\neq 0$ or $h\neq -\frac{1}{2}$, then\[\rk \wti{d}^0_{1,+}|_{h}\ge \frac{1}{2}\left(\dim KHI(-S_0^3(K),\wti{K}_0,h)-\dim KHI(-S_0^3(K),\wti{K}_0,h+ 1)\right).\]If either $\tau_I(K)\neq 0$ or $h\neq \frac{1}{2}$, then \[\rk \wti{d}^0_{1,-}|_{h}\ge \frac{1}{2}\left(\dim KHI(-S_0^3(K),\wti{K}_0,h)-\dim KHI(-S_0^3(K),\wti{K}_0,h- 1)\right).\]
}

\bpf
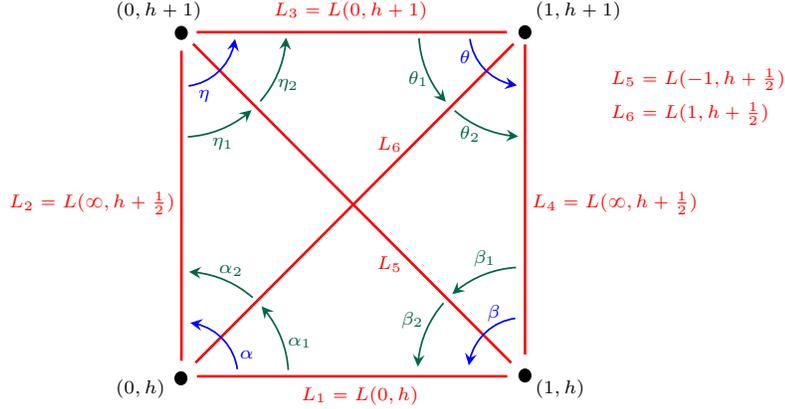
\begin{figure}[ht]
\vspace{-0.2in}
	\centering
	\begin{overpic}[width = 0.6\textwidth]{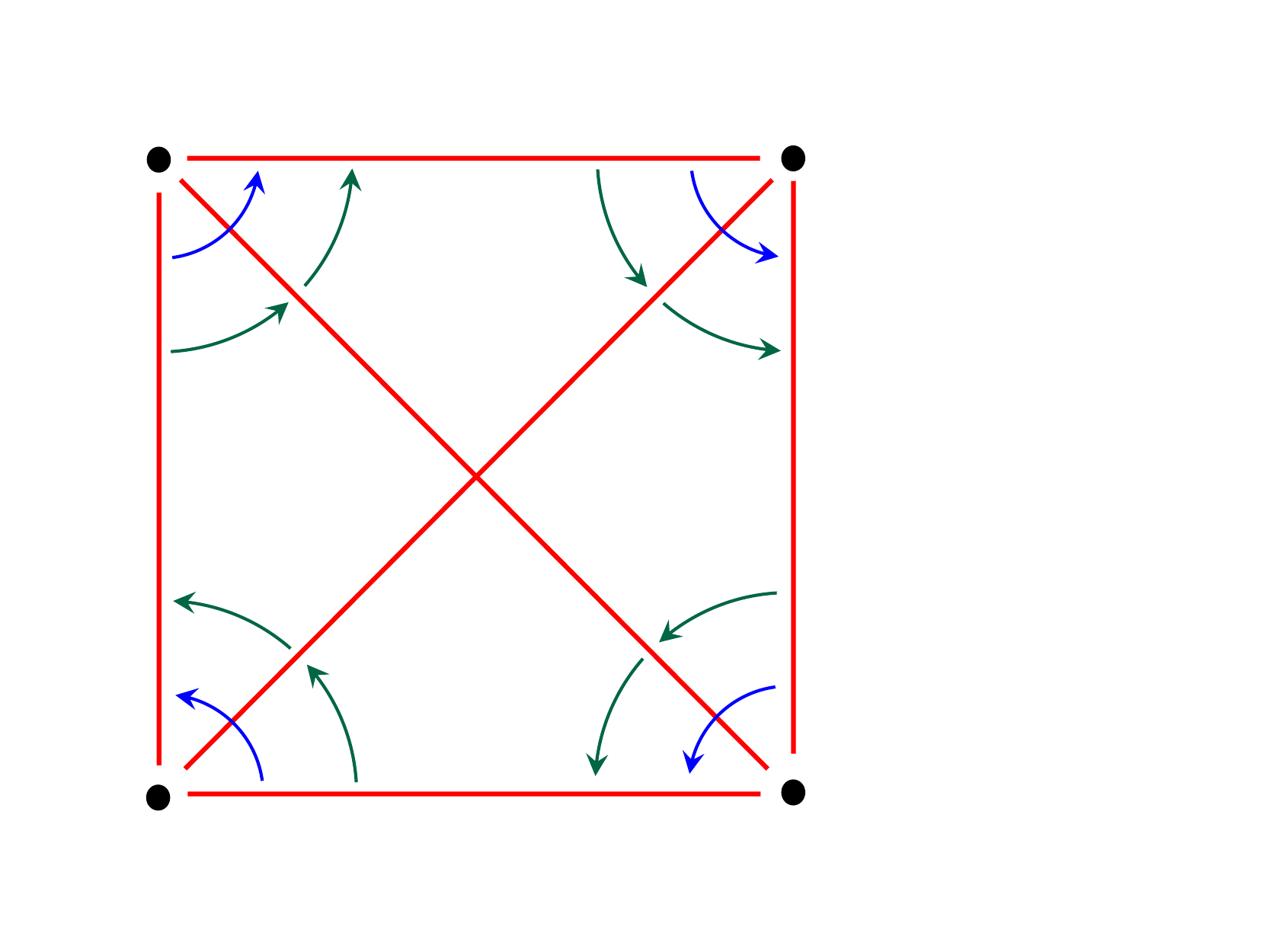}
		\put(3,10) {\tiny $(0,h)$}
		\put(64,10) {\tiny $(1,h)$}
		\put(3,65) {\tiny $(0,h+1)$}
		\put(64,65) {\tiny $(1,h+1)$}
		
		\put(30,9) {\tiny\color{red} $L_1 = L(0,h)$}
		\put(26,65) {\tiny\color{red} $L_3 = L(0,h+1)$}
		\put(-12.5,37) {\tiny\color{red} $L_2 = L(\infty,h+\frac{1}{2})$}
		\put(63.5,37) {\tiny\color{red} $L_4 = L(\infty,h+\frac{1}{2})$}
		\put(75,55) {\tiny\color{red} $L_5 = L(-1,h+\frac{1}{2})$}
		\put(75,50) {\tiny\color{red} $L_6 = L(1,h+\frac{1}{2})$}
		\put(41,45.5) {\tiny\color{red} $L_6$}
		\put(41,28) {\tiny\color{red} $L_5$}
		
		\put(21,15) {\tiny\color{blue} $\alpha$}
		\put(28,17) {\tiny\color{lygreen} $\alpha_1$}
		\put(18,27.5) {\tiny\color{lygreen} $\alpha_2$}
		
		\put(57,21) {\tiny\color{blue} $\beta$}
		\put(55,29) {\tiny\color{lygreen} $\beta_1$}
		\put(44.5,20) {\tiny\color{lygreen} $\beta_2$}
		
		\put(53,58) {\tiny\color{blue} $\theta$}
		\put(45.5,55) {\tiny\color{lygreen} $\theta_1$}
		\put(53,47) {\tiny\color{lygreen} $\theta_2$}
		
		\put(15,53) {\tiny\color{blue} $\eta$}
		\put(17,46) {\tiny\color{lygreen} $\eta_1$}
		\put(26.5,54.5) {\tiny\color{lygreen} $\eta_2$}
		
	\end{overpic}
	\vspace{-0.3in}
	\caption{The graphic illustration for the octahedral diagram.}\label{fig5}
\end{figure}
We will prove the statement only for $\wti{d}^0_{1,+}$. The proof for $\wti{d}^0_{1,-}$ is obtained by taking the orientation reversal of the knot. Following \S \ref{subsec: illustration of bypass maps} \eqref{t1: arc}, we use the notation $L(p/q,h)$ to represent $KHI(-S^3_{p/q},\wti{K}_{p/q},h)$. From \S \ref{subsec: illustration of bypass maps} \eqref{t7: octahedron}, we have the following octahedral diagram, with a graphic illustration in Figure \ref{fig5}.
\begin{equation}\label{eq: octahedron}
	\xymatrix{
	&L_5=L(-1,h+\frac{1}{2})\ar[dr]^{\eta_2}\ar[ddl]^<<<<<{\be_2}&\\
	L_4=L(\infty,h+\frac{1}{2})\ar[ur]^{\be_1}\ar[d]_{\be}&&L_3=L(0,h+1)\ar[ll]^{\theta}\ar[ddl]_>>>>>{\theta_1}\\
	L_1=L(0,h)\ar[dr]_{\al_1}\ar[rr]^{\al}&&L_2=L(\infty,h+\frac{1}{2})\ar[u]_{\eta}\ar[uul]^>>>>>{\eta_1}\\
	&L_6=L(1,h+\frac{1}{2})\ar[uul]_<<<<<{\theta_2}\ar[ur]_{\al_2}&
	}
\end{equation}
where we abuse notation and use the angles to represent the angle maps.
Note that \[\al_1=\psi^{0}_{-,1},~\al_2=\psi^1_{+,\mu},~\al=\psi^0_{+,\mu},~\be_1=\psi^\mu_{-,-1},~\be_2=\psi^{-1}_{+,0},~\be=\psi^\mu_{-,0}\]\[\eta_1=\psi^\mu_{+,-1},~\eta_2=\psi^{-1}_{-,0},\eta=\psi^{\mu}_{+,0},~\theta_1=\psi^{0}_{+,1},~\theta_2=\psi^1_{-,\mu},\aand \theta=\psi^{0}_{-,\mu}.\]

If $\im\be_1\cap \im \eta_1\neq 0$, then Corollary \ref{cor: char of tau G}
implied that $\tau_I(K)=0$ and $h+\frac{1}{2}=0$, which contradicted to the assumption. Then $\im\be_1\cap \im \eta_1=0$. Then from the commutative diagram $\be_1\circ \theta_2=\eta_1\circ \al_2$ (cf. \S \ref{subsec: illustration of bypass maps} \eqref{t7: octahedron} or ultimately \cite[Lemma 4.33]{LY2020}), we know that $\be_1\circ \theta_2=0$. 

From \cite[Lemma 10.3]{KM2016octahedral}, the $4$-periodic complex made by $L_i$ for $i=1,2,3,4$ and the maps $\al,\eta,\theta,\be$ satisfies the condition that the homologies at diametrically opposite corners have the same dimension:
\begin{equation}\label{eq: homology at opposite corner}
	\dim \ker\beta - \rk\theta = \dim \ker\eta - \rk \alpha\text{ and }\dim \ker \theta - \rk\eta = \dim \ker \alpha -\rk\beta.
\end{equation}
Moreover, the proof \cite[Lemma 10.3]{KM2016octahedral} implies that the dimension of the homology at $L_2$ and $L_4$ equals $\rk(\be_1\circ \theta_2)$, which is zero as discussed previously. Hence, we conclude
\[
\dim \ker\beta = \rk \theta\text{ and }\dim \ker \eta = \rk \alpha
\]
and consequently
\[\dim L_2=\rk\al+\rk\eta,~\dim L_4=\rk \be+\rk \theta.\]
Note that by definition $L_2=L_4$, and hence
\[\rk\al+\rk\eta=\rk \be+\rk \Longrightarrow \rk\be-\rk \eta = \rk \al - \rk \theta.\]
Also, from Equation \eqref{eq: homology at opposite corner}, we have
\[
\begin{aligned}
	\dim L_1-\dim L_3 &= (\rk\alpha + \dim \ker\al) - (\rk \theta +\dim \ker\theta)\\
	&=(\rk\alpha + \dim \ker\al - \rk\beta + \rk\beta) - (\rk \theta +\dim \ker\theta - \rk\eta + \rk \eta)\\
	&=\rk \al+\rk\be-\rk \eta-\rk \theta\\
	&= 2(\rk\be-\rk \eta).
\end{aligned}
\]
Since the angle for $\be\circ \al$ is $\pi$, we have $\be\circ \al=\wti{d}^0_{1,+}$ (cf. \S \ref{subsec: illustration of bypass maps} \eqref{t6: differential}). We have a surjective map
\[
\ov{\be}:L_1\slash\im\alpha\to \im\beta\slash\im(\beta\circ\alpha),
\]
which implies that
\[
\begin{aligned}
	\rk \wti{d}^0_{1,+} &= \dim \im(\beta\circ\al)\\
	&=\dim \im\beta - \dim\im\ov{\be}\\
	&=\rk\beta -\left(\dim \left(L_1\slash\im\al\right) - \dim\ker\ov{\be}\right)\\
	&\geq \rk\beta - \dim \left(L_1\slash\im\al\right)\\
	&= \rk\beta - \dim \left(L_1\slash\ker \eta\right)\\
	&=\rk\beta-\rk\eta\\
	&=\frac{1}{2}(\dim L_1-\dim L_3)
\end{aligned}
\]
This completes the proof of the proposition.
\epf
\noindent{\bf Proposition \ref{prop: d0 diff}}
{\it
Suppose $K\subset S^3$ is a knot of genus $g>0$. Then we have \begin{equation*}\label{eq: half rank}
    \rk \wti{d}^0_{1,\pm}|_{\pm (g-\frac{1}{2})}=\frac{1}{2}\dim KHI(-S_0^3(K),\wti{K}_0,\pm (g-\frac{1}{2}))\ge 1.
\end{equation*}
Hence we have
\begin{equation*}\label{eq: homology zero}
    H_*(KHI(-S_0^3(K),\wti{K}_0,\pm (g-\frac{1}{2})),\wti{d}^0_{1,\pm}|_{\pm (g-\frac{1}{2})})=0.
\end{equation*}
}

\bpf
This follows from Proposition \ref{prop: rank d1 at least} and the following facts from \cite[Theorem 2.21(1)(2)]{LY2020}\[KHI(-S_0^3(K),\wti{K}_0,\pm (g+\frac{1}{2}))=0\aand KHI(-S_0^3(K),\wti{K}_0,\pm (g-\frac{1}{2}))\neq 0.\] 
\epf

\subsection{The main theorem}
In this subsection, we prove Theorem \ref{thm: main}. Recall that for a closed, oriented $3$-manifold $Y$, we define\begin{equation}\label{eq: defn of t_2}
	t_2(Y) = \frac{1}{2}\left(\dim I^{\sharp}(Y;\F_2=\Z/2) - \dim I^{\sharp}(Y;\C)\right).
\end{equation}

A key observation is the following.
\blem\label{lem: connected sum with lens space}
	Suppose $Y$ is a closed oriented $3$-manifold and $L(p,q)$ is a lens space. If $t_2(Y)=0$, then $t_2(Y\# L(p,q))=0$.
\elem
\bpf
From Remark \ref{rem: mirror knot trick}, we know that $t_2(Y)=0$ if and only if $t_2(-Y)=0$. Since \[-\left(Y\# L(p,q)\right)=(-Y)\#\left(-L(p,q)\right)=(-Y)\#L(p,-q),\]it suffices to prove the lemma for $p,q>0$. 

From now on, we assume $p,q>0$. From \cite[Corollary 1.2]{Scaduto2015instanton} and the fact that a lens space is a double cover of a rational knot (that is a quasi-alternating knot), we have \begin{equation*}\label{eq: rk of L(p,q)}
    \dim I^\sharp(L(p,q);\K)=p
\end{equation*}for any field $\K$, especially for $\K=\F_2$ or $\C$. From the connected sum formula in \cite[Remark 1.6]{li2018contact}, we have\[
\dim I^{\sharp}(Y\# L(p,q);\C) =\dim I^{\sharp}(Y;\C)\cdot \dim I^\sharp(L(p,q);\C)= p\cdot I^{\sharp}(Y;\C).
\]

Following the procedure in the proofs of \cite[Proposition 7.3]{kronheimer2007monopolesandlens} and \cite[Proposition 4.3]{baldwin2020concordance}, for fixed pair $(p,q)$ with $p,q>0$, one can find two other pairs of co-prime non-negative integers $(p^\p,q^\p)$ and $(p^\pp,q^\pp)$ such that
\[
p = p^\p + p^\pp \aand q = q^\p + q^\pp,
\]
and they form a surgery triad (along an unknot in $S^3$):
\[
\xymatrix{
L(p^\p,q^\p)\ar[rr]&&L(p^\pp,q^\pp)\ar[dl]\\
&L(p,q)\ar[lu]&
}
\]
Since the surgery operation is local, this triad is preserved when taking the connected sum with $Y$. From \cite[\S 7.5]{Scaduto2015instanton}, we obtain an exact triangle for framed instanton homology
\[
\xymatrix{
I^{\sharp}(Y\#L(p^\p,q^\p);\F_2)\ar[rr]&&I^{\sharp}(Y\#L(p^\pp,q^\pp);\F_2)\ar[dl]\\
&I^{\sharp}(Y\#L(p,q);\F_2)\ar[lu]&
}
\]
We prove the lemma by induction: the fact that $L(1,q)\cong S^3$ for any $q$ serves as the initial point, and if $t_2=0$ for the connected sum already holds for $L(p^\p,q^\p)$ and $L(p^\pp,q^\pp)$, then
\[
\begin{aligned}
	\dim I^{\sharp}(Y\#L(p,q);\C) &\leq \dim I^{\sharp}(Y\#L(p,q);\F_2)\\
	&\leq \dim I^{\sharp}(Y\#L(p^\p,q^\p);\F_2) + \dim I^{\sharp}(Y\#L(p^\pp,q^\pp);\F_2)\\
	&= \dim I^{\sharp}(Y\#L(p^\p,q^\p);\C) + \dim I^{\sharp}(Y\#L(p^\pp,q^\pp);\C)\\
	&= (p^\p+p^\pp)\cdot \dim I^{\sharp}(Y;\C)\cdot \\
	&=p\cdot \dim I^{\sharp}(Y;\C)\\
	&= \dim I^{\sharp}(Y\#L(p,q);\C).
\end{aligned}
\]
Thus, we conclude that $t_2(Y\#L(p,q);\F_2) = 0$.
\epf

\blem\label{lem: t_2 preserves until nu-sharp}
Suppose $K\subset S^3$ is a knot. For any $n\in\Z$, we have
\begin{equation}\label{eq: t_2 drops by at most 1}
	t_2(S^3_{n+1}(K)) \geq t_2(S^3_{n}(K)) - 1.
\end{equation}
Furthermore, if $\nu^{\sharp}(K)>0$ and $n\in [1,\nu^{\sharp}(K)]$, then we have
\[
t_2(S^3_n(K))\geq t_2(S^3_1(K)).
\]
\elem
\bpf
From \cite[\S 7.5]{Scaduto2015instanton}, for any field $\K$ and any $n\in\Z$, there exists a surgery exact triangle
\[
\xymatrix{
I^{\sharp}(S^3_n(K);\K)\ar[rr]&&I^{\sharp}(S^3_{n+1}(K);\K)\ar[dl]\\
&I^{\sharp}(S^3;\K)\ar[lu]&
}
\]
Since $\dim I^{\sharp}(S^3;\K) = 1$, we know that
\begin{equation}\label{eq: dim differed by 1}
	\dim I^{\sharp}(S^3_{n+1}(K);\K) = \dim I^{\sharp}(S^3_{n}(K);\K) \pm 1.
\end{equation}
Comparing this with the definition of $t_2$ in \eqref{eq: defn of t_2}, we conclude that
\[
t_2(S^3_{n+1}(K)) \geq t_2(S^3_{n}(K)) - 1.
\]

From \cite[Theorem 1.1]{baldwin2020concordance}, for any $n\in [1,\nu^{\sharp}-1]$, we have
\[
\dim I^{\sharp}(S^3_{n+1}(K);\C) = \dim I^{\sharp}(S^3_{n}(K);\C) - 1.
\]
Thus, we conclude from \eqref{eq: dim differed by 1} and \eqref{eq: defn of t_2} that
\[
t_2(S^3_{n+1}(K))\geq t_2(S^3_n(K))\geq\cdots\geq t_2(S^3_1(K)).
\] 
\epf

\blem\label{lem: dim of 2g-1 dual knot is strictly larger}
Suppose $K\subset S^3$ is an instanton L-space knot of genus $g$. Then 
\begin{equation}\label{eq: KHI 2g-1}
    \dim KHI(S^3_{2g-1}(K),\wti{K}_{2g-1})=2g+1.
\end{equation}
and $\tau_I(K)=g$.
\elem
\bpf
From \cite[Theorem 1.15]{baldwin2019lspace}, \cite[Theorem 1.17 and Remark 1.19]{LY2021large}, and \cite[Proposition 1.5]{LY2020}, for an instanton L-space knot $K$, when $n\in\Z$ and $n>2g(K)-1$, we have\[\dim KHI(S^3_{n}(K),\wti{K}_{n})=\dim I^\sharp(S^3_n(K);\C)=n\]\[\aand \dim KHI(S^3_{2g-1}(K),\wti{K}_{2g-1})>\dim I^\sharp(S^3_{2g-1}(K);\C)=2g-1.\]From Lemma \ref{lem: V-shaped}, we obtain \eqref{eq: KHI 2g-1} and $\tau_I(K)=g$.

\epf

Now, it is time to prove the main theorem.

\noindent{\bf Theorem \ref{thm: main}}
{\it
	Suppose $K\subset S^3$ is a non-trivial knot and $r\in\Q_+$ such that $I^{\sharp}(S^3_r(K);\Z)$ has no $2$-torsion. Then $K$ is an instanton L-space knot, and $r> 2g(K) - 1$. Furthermore, in the following special cases, we have a better lower bound:
	\begin{itemize}
		\item If $r\in\Z_+$, then we have
		\[
			r\geq 2g(K) - 1 + t_2(S^3_1(K)),
		\]
		where for any closed oriented $3$-manifold $Y$, we define
		\begin{equation*}
		    t_2(Y)=\frac{1}{2}\left(\dim I^\sharp(Y;\F_2=\Z/2)-\dim I^\sharp(Y;\C)\right).
		\end{equation*}Note that $t_2(S^3_1(K))\ge 1$ by \cite[Theorem 1.1]{LY2025torsion}.
		\item If $r\in\Z_+-\frac{1}{2}$, then we have \[r> 2g(K).\]
	\end{itemize}
}
\bpf
Suppose we have a non-trivial knot $K\subset S^3$ of genus $g$ and $r\in\Q_+$ such that $t_2(S^3_r(K)) = 0$. 

We first assume that $r=n\in \Z_+$. Then Lemma \ref{lem: 2-tor to C} and Proposition \ref{prop: 2-torsion free implies being L-space} imply that $K$ is an instanton L-space knot and $n+1>2g-1$. From \cite[Theorem 1.16]{baldwin2020concordance}, we have $\nu^{\sharp}(K) = 2g - 1$. Lemma \ref{lem: t_2 preserves until nu-sharp} and the assumption $t_2(S^3_n(K)) = 0$ then imply
\[
n - (2g-1)\geq t_2(S^3_{2g-1}(K))-t_2(S^3_n(K))=t_2(S^3_{2g-1}(K))  \geq t_2(S^3_1(K)),
\]
which implies
\[
n\geq 2g - 1 + t_2(S^3_1(K)).
\]From \cite[Theorem 1.1]{LY2025torsion}, we have $t_2(S^3_1(K))\ge 1$. Hence $n> 2g-1$.

Then we assume that $r=p/q$ with $p>0$ and $q>1$. From \eqref{eq: diffeo for cabling} and Lemma \ref{lem: connected sum with lens space}, we have $t_2(S^3_{pq}(K)) = 0$. Then Lemma \ref{lem: 2-tor to C} and Proposition \ref{prop: 2-torsion free implies being L-space} imply that $K_{p,q}$ is an instanton L-space knot.
Then \cite[Theorem 1.18]{baldwin2020concordance} implies that $K$ is an instanton L-space knot and $p/q>2g(K)-1$. 

Finally, we assume that $r\in \Z_+-\frac{1}{2}$. The above discussion already showed that $r>2g-1$. Thus, it suffices to show that \[t_2\left(S^3_{\frac{4g-1}{2}}(K)\right)\neq 0.\]We fix the meridian $\mu$ and the Seifert framing $\lambda$ of $K$ with $\mu\cdot \lambda =-1$. A framing of the dual knot $\wti{K}_{2g-1}\subset S^3_{2g-1}(K)$ can then be expressed as
\[
\wti{\mu} = (2g-1)\mu + \lambda \aand \wti{\lambda} = -(4g-1)\mu - 2 \lambda.
\]
Then
\[
\wti{\lambda} + 2\wti{\mu} = \mu.
\]
Hence we obtain by Bhat's exact triangle \cite[Theorem 1.1]{bhat2023newtriangle} that
\[
\xymatrix{
I^{\sharp}(S^3_{\frac{4g-1}{2}}(K);\F_2)\ar[rr] && I^{\sharp}(S^3;\F_2)\ar[dl]\\
&I^{\sharp}(S^3_{2g-1}(K),\wti{K}_{2g-1};\F_2)\ar[lu]&
}
\]
Thus, together with \cite[Formula (1.3)]{LY2025torsion}, the fact that $K$ is an instanton L-space knot, Lemma \ref{lem: dim of 2g-1 dual knot is strictly larger}, and \cite[Theorem 1.16]{baldwin2020concordance}, we have
\[
\begin{aligned}
\dim I^{\sharp}(S^3_{\frac{4g-1}{2}}(K);\F_2) &\geq \dim I^{\sharp}(S^3_{2g-1}(K),\wti{K};\F_2) - 1
\\&\geq 2\cdot \dim  KHI(S^3_{2g-1}(K),\wti{K}_{2g-1}) - 1\\
&= 4g + 1\\
&= \dim  I^{\sharp}\left(S^3_{\frac{4g-1}{2}}(K);\C\right) + 2.
\end{aligned}
\]
Hence we conclude that
\[
t_2\left(S^3_{\frac{4g-1}{2}}(K)\right)\neq 0.
\]
\epf

\subsection{Corollaries}\label{subsec: corollaries}
In this subsection, we prove corollaries in \S \ref{sec: introduction}.

\noindent{\bf Corollary \ref{cor: small surgery}}
{\it
	Suppose $K\subset S^3$ is a non-trivial knot, and $r\in \Q$ with $|r|\in (0,2g-1]$. Then $I^{\sharp}(S^3_{r}(K);\Z)$ admits $2$-torsion.
}
\bpf
This follows directly from Theorem \ref{thm: main} for $r>0$ and Remark \ref{rem: mirror knot trick} for $r<0$.
\epf

\noindent{\bf Corollary \ref{cor: 2g-1}}
{\it
	Suppose $K\subset S^3$ is a non-trivial knot. Suppose $r=p/q$ with $p,q>0$ and $\Delta_K(\zeta^2) \neq 0$ for any $p$-th root of unity $\zeta$, where $\Delta_K(t)$ is the Alexander polynomial of $K$ (indeed, when $p$ is a prime power, the condition for roots always holds by \cite[Remark 1.6]{baldwin2016contact}). Then $S^3_r(K)$ is non-$SU(2)$-abelian unless
	\[
	r > 2g(K) - 1.
	\]
	Furthermore, if $r\in\Z_+$, then
	\[
	r \geq 2g(K) - 1 + t_2(S^3_1(K))
	\]
	and if $r\in\Z_+-\frac{1}{2}$, then
	\[
	r> 2g(K).
	\]
}

\bpf
The proof is similar to \cite[Theorem 1.5]{baldwin2019lspace}. Assume that $r=p/q$ with $p,q>0$ such that $Y=S^3_r(K)$ is $SU(2)$-abelian. Since any $SU(2)$-representation has an abelian image, it factors through $H_1(Y;\Z)$. Thus, there are precisely $p=|H_1(Y)|$ such representations. All of them are isolated, and the assumption on the $p$-th root of unity guarantees the required non-degeneracy condition. Thus one can construct an instanton chain complex of $I^{\sharp}(S^3_r(K))$ generated by precisely those $p$ reducible representations. Since Scaduto \cite[Corollary 1.4]{Scaduto2015instanton} showed that the Euler characteristic $\chi(I^{\sharp}(S^3_r(K)))=p$, we know there is no differential over any coefficients. Thus,
\[
t_2(S^3_r(K)) = 0.
\]
Then Theorem \ref{thm: main} applies and we obtain the corollary.
\epf

\noindent{\bf Corollary \ref{cor: 5- and 5.5-surgery}}
{\it
	Suppose $K\subset S^3$ is a knot and $r\in\frac{1}{2}\Z$ with $|r|<6$ such that $S^3_r(K)$ is $SU(2)$-abelian. Then one of the following cases happens.
	\begin{itemize}
		\item $K$ is the unknot, 
		\item $K$ is the right-handed trefoil and $r\in\{5,\frac{11}{2}\}$.
            \item $K$ is the left-handed trefoil and $r\in \{-5,-\frac{11}{2}\}$.
	\end{itemize}
}

\bpf
Because of Theorem \ref{thm: BLSY results} and Remark \ref{rem: mirror knot trick}, it suffices to consider $r\in \{5,\frac{11}{2}\}$.

First, assume that $S^3_5(K)$ is $SU(2)$-abelian. Since $5$ is a prime number, the hypothesis of Corollary \ref{cor: 2g-1} is satisfied, which implies that $5>2g(K)-1$, and consequently that $g(K)\leq 2$. Moreover, from the proof of Corollary \ref{cor: 2g-1}, we know that $K$ is an instanton L-space knot. 

From \cite[Theorem 1.5 and Remark 1.7]{baldwin2016contact}, when $g(K)\leq 1$, the unknot and the right-handed trefoil are the only instanton L-space knots. Since their $5$-surgeries are both lens spaces, they are indeed candidates for $K$.

If $g(K) = 2$, we apply \cite[Corollary 1.4]{farber2024fixed} and conclude that $K=T_{2,5}$ is the only candidate. Yet by \cite[Proposition 3.1]{moser1971elementary}, for $r=5$, $s=2$, $p=1$, and $q=-5$, we have that $S^3_5(T_{2,5})\cong \Sigma(2,5,5)$ is a Brieskorn sphere, which, by \cite[Theorem 1.2]{sivek2022menagerie}, is not $SU(2)$-abelian. This concludes the proof of the case $r=5$.

The proof of the case $r= 11/2$ is similar. Since $11$ is a prime number, the hypothesis of Corollary \ref{cor: 2g-1} is again fulfilled. Thus, Corollary \ref{cor: 2g-1} implies that $11/2 \geq 2g(K)$, which further implies that $g(K) \leq 2$. As above, the only three candidates are the unknot, the right-handed trefoil, and the torus knot $T_{2,5}$. For the unknot, the $11/2$-surgery yields the lens space $L(11,-2)$. For the right-handed trefoil, \cite[Proposition 3.1]{moser1971elementary} for $r=3,s=2,p=2,q=-11$ implies that $11/2$-surgery yields the lens space $L(11,8)$. Those two lens spaces are both $SU(2)$-abelian. For $T_{2,5}$, \cite[Proposition 3.1]{moser1971elementary} for $r=5,s=2,p=2,q=-11$ implies $S^3_{11/2}(T_{2,5})\cong \Sigma (2,5,9)$, which again by \cite[Theorem 1.2]{sivek2022menagerie} is not $SU(2)$-abelian.
\epf

\noindent{\bf Corollary \ref{cor: 7-surgery}}
{\it
		Suppose $K\subset S^3$ is a knot such that $S^3_7(K)$ is $SU(2)$-abelian, then either
	\begin{itemize}
		\item $K$ is the unknot or the right-handed trefoil, or
		\item $K$ is a genus-$3$ instanton L-space knot with $t_2(S^3_1(K))\leq 2$.
	\end{itemize}
}
\bpf
The proof is similar to that of Corollary \ref{cor: 5- and 5.5-surgery}. If $7$-surgery on a knot $K\subset S^3$ is $SU(2)$-abelian, then $K$ is a instanton L-space knot and
\[
2g(K) - 1 + t_2(S^3_1(K)) \leq 7.
\]
Hence, either $g(K) = 3$, and $t_2(S^3_1(K))\leq 2$, or $g(K)\leq 2$, so $K$ is among the unknot, the right-handed trefoil, and the torus knot $T_{2,5}$. We can then similarly exclude the last candidate by applying \cite[Proposition 3.1]{moser1971elementary} for $r=5,s=2,p=1,q=-7$ to obtain $S^3_{7}(T_{2,5})\cong \Sigma (2,5,3)$, and then applying \cite[Theorem 1.2]{sivek2022menagerie}. Note that \cite[Proposition 3.2]{moser1971elementary} for $r=3,s=2,p=1,q=-7$ implies that $S_7^3(T{2,3})\cong L(7,4)$, which is $SU(2)$-abelian.
\epf
\section{Future direction}\label{sec: future direction}
In this section, we provide some conjectures and open questions for future research. First, \cite[Theorem 1.5]{bhat2023newtriangle} (and also Remark \ref{rem: small t_2}) showed that $I^\sharp(P;\Z)$ has $2$-torsion for the Poincar\'e sphere $P$. Thus, we make the following conjecture, as an instanton analogue of the Poincar\'e conjecture.
\begin{conjecture}\label{conj: su2}
	The $3$-sphere $S^3$ is the only connected, closed, oriented $3$-manifold $Y$ with $\dim I^\sharp(Y;\F_2)=1$.
\end{conjecture}

Note that our main theorem, Theorem \ref{thm: main}, together with \[\chi(I^\sharp(S^3_{p/q}(K);\F_2))=|p|\]by \cite[Corollary 1.4]{Scaduto2015instanton}, indicate that Conjecture \ref{conj: su2} holds for all $3$-manifolds coming from surgery on knots in $S^3$. Also note that Conjecture \ref{conj: su2} implies the classical $3$-dimensional Poincar\'e conjecture.

Another thing to note is that currently, all examples for $t_2(Y) = 0$ comes from instanton L-spaces $Y$. So we would like to propose the following question.
\begin{question}
	Can we find a closed oriented $3$-manifold $Y$ such that $t_2(Y) = 0$ and $Y$ is not an instanton-L-space?
\end{question}

Another question is from Remark \ref{rem: small t_2}.
\begin{question}
	Can we find a knot $K\subset S^3$ with $t_2(S^3_1(K))\in \{1,2\}$?
\end{question}

At last, observe that the combination of results in \cite{daemi2022ribbon,gordon1981ribbon,baldwin2019lspace} imply that instanton L-space knots are minimal under the partial ordering induced by ribbon concordance. Thus, our main theorem, Theorem \ref{thm: main}, directly relates the fact of admitting $2$-torsion-free Dehn surgeries with this minimality. Hence, it is natural to ask the following question.
\begin{question}
	Can we use the $2$-torsion technique to define new effective obstructions to ribbon concordance between knots and (rational) ribbon cobordisms between connected, closed, oriented $3$-manifolds?
\end{question}
\bibliographystyle{alpha}

\begin{thebibliography}{DLVVW22}

\bibitem[Bha24]{bhat2023newtriangle}
Deeparaj Bhat.
\newblock Surgery exact triangles in instanton theory.
\newblock ArXiv: 2311.04242, v2, 2024.

\bibitem[BLSY24]{BLSY21}
John~A. Baldwin, Zhenkun Li, Steven Sivek, and Fan Ye.
\newblock Small {D}ehn surgery and {${\rm SU}(2)$}.
\newblock {\em Geom. Topol.}, 28(4):1891--1922, 2024.

\bibitem[BS15]{BS2015naturality}
John~A. Baldwin and Steven Sivek.
\newblock Naturality in sutured monopole and instanton homology.
\newblock {\em J. Differential Geom.}, 100(3):395--480, 2015.

\bibitem[BS16a]{baldwin2016contact}
John~A. Baldwin and Steven Sivek.
\newblock A contact invariant in sutured monopole homology.
\newblock {\em Forum Math. Sigma}, 4:e12, 82, 2016.

\bibitem[BS16b]{BS2016instanton}
John~A. Baldwin and Steven Sivek.
\newblock Instanton {F}loer homology and contact structures.
\newblock {\em Selecta Math. (N.S.)}, 22(2):939--978, 2016.

\bibitem[BS21]{baldwin2020concordance}
John~A. Baldwin and Steven Sivek.
\newblock Framed instanton homology and concordance.
\newblock {\em J. Topol.}, 14(4):1113--1175, 2021.

\bibitem[BS22]{baldwin2022concordanceII}
John~A. Baldwin and Steven Sivek.
\newblock {Framed instanton homology and concordance II}.
\newblock arXiv: 2206.11531, v1, 2022.

\bibitem[BS23]{baldwin2019lspace}
John~A. Baldwin and Steven Sivek.
\newblock Instantons and {L}-space surgeries.
\newblock {\em J. Eur. Math. Soc. (JEMS)}, 25(10):4033--4122, 2023.

\bibitem[DLVVW22]{daemi2022ribbon}
Aliakbar Daemi, Tye Lidman, David~Shea Vela-Vick, and C.-M.~Michael Wong.
\newblock Ribbon homology cobordisms.
\newblock {\em Adv. Math.}, 408:Paper No. 108580, 68, 2022.

\bibitem[FRW24]{farber2024fixed}
Ethan Farber, Braeden Reinoso, and Luya Wang.
\newblock Fixed-point-free pseudo-{A}nosov homeomorphisms, knot {F}loer homology and the cinquefoil.
\newblock {\em Geom. Topol.}, 28(9):4337--4381, 2024.

\bibitem[GL23]{li2019decomposition}
Sudipta Ghosh and Zhenkun Li.
\newblock Decomposing sutured monopole and instanton {F}loer homologies.
\newblock {\em Selecta Math. (N.S.)}, 29(3):Paper No. 40, 60, 2023.

\bibitem[GLW24]{ghosh2024tau}
Sudipta Ghosh, Zhenkun Li, and C.-M.~Michael Wong.
\newblock On the tau invariants in instanton and monopole {F}loer theories.
\newblock {\em J. Topol.}, 17(2):Paper No. e12346, 53, 2024.

\bibitem[Gor81]{gordon1981ribbon}
C.~McA. Gordon.
\newblock Ribbon concordance of knots in the {$3$}-sphere.
\newblock {\em Math. Ann.}, 257(2):157--170, 1981.

\bibitem[Gor83]{gordon1983dehn}
C.~McA. Gordon.
\newblock Dehn surgery and satellite knots.
\newblock {\em Trans. Amer. Math. Soc.}, 275(2):687--708, 1983.

\bibitem[HKM09]{honda2009contactInSFH}
Ko~Honda, William~H. Kazez, and Gordana Mati\'c.
\newblock The contact invariant in sutured {F}loer homology.
\newblock {\em Invent. Math.}, 176(3):637--676, 2009.

\bibitem[Hom14]{Hom2014epsilon}
Jennifer Hom.
\newblock The knot {F}loer complex and the smooth concordance group.
\newblock {\em Comment. Math. Helv.}, 89(3):537--570, 2014.

\bibitem[Hon00a]{honda2000classification}
Ko~Honda.
\newblock On the classification of tight contact structures {I}.
\newblock {\em Geom. Topol.}, 4:309--368, 2000.

\bibitem[Hon00b]{Honda2000classification2}
Ko~Honda.
\newblock On the classification of tight contact structures. {II}.
\newblock {\em J. Differential Geom.}, 55(1):83--143, 2000.

\bibitem[Juh06]{juhasz2006holomorphic}
Andr\'{a}s Juh\'{a}sz.
\newblock Holomorphic discs and sutured manifolds.
\newblock {\em Algebr. Geom. Topol.}, 6:1429--1457, 2006.

\bibitem[Kir97]{kirby1997problems}
Rob Kirby.
\newblock Problems in low-dimensional topology.
\newblock In {\em Geometric topology ({A}thens, {GA}, 1993)}, volume 2.2 of {\em AMS/IP Stud. Adv. Math.}, pages 35--473. Amer. Math. Soc., Providence, RI, 1997.

\bibitem[KM04]{kronheimer04su2}
Peter~B. Kronheimer and Tomasz~S. Mrowka.
\newblock {Dehn surgery, the fundamental group and {SU}(2)}.
\newblock {\em Math. Res. Lett.}, 11(5-6):741--754, 2004.

\bibitem[KM10]{kronheimer2010knots}
Peter~B. Kronheimer and Tomasz~S. Mrowka.
\newblock Knots, sutures, and excision.
\newblock {\em J. Differ. Geom.}, 84(2):301--364, 2010.

\bibitem[KM16]{KM2016octahedral}
P.~B. Kronheimer and T.~S. Mrowka.
\newblock Exact triangles for {$SO(3)$} instanton homology of webs.
\newblock {\em J. Topol.}, 9(3):774--796, 2016.

\bibitem[KMOS07]{kronheimer2007monopolesandlens}
Peter~B. Kronheimer, Tomasz~S. Mrowka, Peter~S. Ozsv\'{a}th, and Zolt\'{a}n Szab\'{o}.
\newblock Monopoles and lens space surgeries.
\newblock {\em Ann. of Math. (2)}, 165(2):457--546, 2007.

\bibitem[Li20]{li2018contact}
Zhenkun Li.
\newblock Contact structures, excisions and sutured monopole {F}loer homology.
\newblock {\em Algebr. Geom. Topol.}, 20(5):2553--2588, 2020.

\bibitem[Li21a]{li2018gluing}
Zhenkun Li.
\newblock Gluing maps and cobordism maps in sutured monopole and instanton {F}loer theories.
\newblock {\em Algebr. Geom. Topol.}, 21(6):3019--3071, 2021.

\bibitem[Li21b]{li2019direct}
Zhenkun Li.
\newblock Knot homologies in monopole and instanton theories via sutures.
\newblock {\em J. Symplectic Geom.}, 19(6):1339--1420, 2021.

\bibitem[LY21]{LY2021large}
Zhenkun Li and Fan Ye.
\newblock ${SU}$$(2)$ representations and a large surgery formula.
\newblock arXiv:2107.11005, v1, 2021.

\bibitem[LY22]{LY2020}
Zhenkun Li and Fan Ye.
\newblock Instanton {F}loer homology, sutures, and {H}eegaard diagrams.
\newblock {\em J. Topol.}, 15(1):39--107, 2022.

\bibitem[LY24]{LY2022integral1}
Zhenkun Li and Fan Ye.
\newblock {Knot surgery formulae for instanton Floer homology I: the main theorem}.
\newblock arXiv:2206.10077, v3, 2024.

\bibitem[LY25a]{LY2025torsion}
Zhenkun Li and Fan Ye.
\newblock 2-torsion in instanton {F}loer homology.
\newblock {\em Adv. Math.}, 472:Paper No. 110289, 55, 2025.

\bibitem[LY25b]{LY2022integral2}
Zhenkun Li and Fan Ye.
\newblock Knot surgery formulae for instanton {F}loer homology {II}: applications.
\newblock {\em Math. Ann.}, 391(4):6291--6371, 2025.

\bibitem[Mos71]{moser1971elementary}
Louise Moser.
\newblock Elementary surgery along a torus knot.
\newblock {\em Pacific J. Math.}, 38:737--745, 1971.

\bibitem[Sar15]{sarkar15moving}
Sucharit Sarkar.
\newblock {Moving basepoints and the induced automorphisms of link Floer homology}.
\newblock {\em Algebr. Geom. Topol.}, 15(5):2479--2515, 2015.

\bibitem[Sca15]{Scaduto2015instanton}
Christopher~W. Scaduto.
\newblock Instantons and odd {K}hovanov homology.
\newblock {\em J. Topol.}, 8(3):744--810, 2015.

\bibitem[SZ22]{sivek2022menagerie}
Steven Sivek and Raphael Zentner.
\newblock A menagerie of {$SU(2)$}-cyclic 3-manifolds.
\newblock {\em Int. Math. Res. Not. IMRN}, 2022(11):8038--8085, 2022.

\end{thebibliography}

\end{document}